\documentclass[10pt,journal]{IEEEtran}
\usepackage[T1]{fontenc}
\usepackage[utf8]{inputenc}
\usepackage{amsmath,amssymb,amsthm,mathtools}
\usepackage{cite,booktabs,graphicx,microtype}
\usepackage{placeins}
\usepackage[hidelinks]{hyperref}
\hypersetup{pdftitle={Dual-Mode Extremum Seeking for a Class of Unknown Controlled Diffusions},pdfauthor={Martin Guay, Shimin Wang, Edward Z. Cao, Stefan Damjanovic}}
\newtheorem{theorem}{Theorem}
\newtheorem{proposition}{Proposition}
\newtheorem{lemma}{Lemma}
\newtheorem{corollary}{Corollary}
\theoremstyle{definition}
\newtheorem{assumption}{Assumption}
\newtheorem{definition}{Definition}
\theoremstyle{remark}
\newtheorem{remark}{Remark}
\newcommand{\R}{\mathbb R}
\newcommand{\E}{\mathbb E}
\newcommand{\Prob}{\mathbb P}
\newcommand{\dd}{\mathrm d}
\newcommand{\eps}{\varepsilon}
\newcommand{\transpose}{^{\mathsf T}}
\newcommand{\Torus}{\mathbb T}
\newcommand{\Spp}{\mathbb S_{++}}
\newcommand{\cD}{\mathcal D}
\newcommand{\cB}{\mathcal B}

\newcommand{\cL}{\mathcal L}
\newcommand{\cA}{\mathcal A}

\newcommand{\ind}{\mathbf 1}
\newcommand{\avg}[1]{\left\langle #1\right\rangle}

\DeclareMathOperator{\col}{col}
\DeclareMathOperator{\diag}{diag}
\DeclareMathOperator{\sym}{sym}
\DeclareMathOperator{\tr}{tr}

\ifCLASSOPTIONtwocolumn

\else

\fi
\allowdisplaybreaks[1]
\newcommand{\phiT}{\varphi}

\begin{document}
\title{Dual-Mode Extremum Seeking for a Class of Unknown Controlled Diffusions}
\author{Martin Guay, Shimin Wang, Edward Z. Cao, Stefan Damjanovic\thanks{Martin Guay and Edward Z. Cao are with the Department of Chemical Engineering, Queen's University, Kingston, ON, Canada. Shimin Wang is with School of Data Science, Lingnan University, Hong Kong. Stefan Damjanovic is with the Department of Electrical and Computer Engineering, University of Toronto, Toronto, ON, Canada}}
\maketitle
\begin{abstract}
A dual-mode extremum-seeking controller is developed for unknown nonlinear
systems with affine control inputs driven by persistent Brownian process
noise. The controller uses a single measured cost, without state estimation or
differentiation of the measurement. Bounded cosine feedback generates, through
Lie-bracket averaging, a proportional--integral (PI) action: the proportional
action stabilizes a feasible minimum of the cost on the time scale of the
plant, and the integral action adapts a bias toward the unknown input that
regulates the plant at that minimum. The averaged closed loop is a PI-regulated diffusion that retains the process noise, and its stability
condition does not require knowledge of the control direction. Under a
regional Lyapunov condition on the averaged dynamics, an exact periodic
response of the noise-free closed loop is constructed and serves as the center
of a Lyapunov function with bounds uniform in the excitation frequency. The
pre-exit mean-square regulation error then consists of an exponentially decaying
transient and a term proportional to the noise variance plus the inverse
excitation frequency, separating noise-induced deviations from probing
oscillations. Exit-probability estimates show that, for small noise and large
excitation frequency, trajectories reach a prescribed neighborhood of the
optimum and remain there for times exponential in the inverse noise variance,
with exponentially small failure probability. Under an additional stochastic
contraction condition, the closed loop is approximated uniformly in time, in
mean square before common exit, by its averaged diffusion driven by the same Brownian motion,
with an error that decreases with the excitation frequency. The regional regulation conditions are verified for a planar hydrodynamic trap
with an unknown stagnation-point offset and actuator sign. Simulations
illustrate intensity-only confinement of a Brownian particle, integral
compensation of the holding input, and axial-strip exit statistics.
\end{abstract}
\section{Introduction}
\label{secI}

Extremum seeking control (ESC) steers an unknown dynamical system toward an
extremum of a measured objective function whose analytical form is unavailable
\cite{AriyurKrstic2003,Scheinker2024}. In its classical form, the plant is
treated as a stable steady-state input--output map that is optimized on a slow
time scale, and closed-loop stability follows from averaging and singular
perturbation arguments \cite{KrsticWang2000,Tan2006}. The time-scale separation
underlying this analysis limits the achievable transient performance and
excludes plants that are not stable under constant inputs. Dual-mode ESC was
introduced to remove this restriction \cite{Guay2016,GuayDochain2017}. It is a
proportional--integral (PI) formulation: a proportional action, generated from
the measured cost, stabilizes the plant at a feasible optimum on the time scale
of the plant, while an integral action adapts a bias toward the unknown constant
input that sustains the optimum at equilibrium. In \cite{GuayAtta2018}, both
actions are realized by periodic bounded feedback and the closed loop is
analyzed through Lie-bracket averaging, which yields practical asymptotic
stability of the unknown optimum for deterministic nonlinear plants. Finite-time
and prescribed-time dual-mode designs are developed in
\cite{Guay2021,Guay2023,Waterman2024}.

Lie-bracket approximation characterizes the effective vector fields generated
by high-frequency, large-amplitude periodic inputs \cite{Durr2013} and has
become a standard framework for the analysis and design of extremum seeking
systems
\cite{Durr2017,Grushkovskaya2018,Grushkovskaya2021,AbdelgalilPoveda2025,Pokhrel2026}.
A bounded trigonometric dependence on the measured cost
\cite{Scheinker2014,Scheinker2016} keeps the update rate bounded irrespective of
the cost magnitude, and the same mechanism provides model-free stabilizing
feedback under unknown control directions through minimum seeking of control
Lyapunov functions \cite{ScheinkerKrstic2013,ScheinkerKrstic2017}. The sign
independence of Lie-bracket feedback has since been exploited to solve robust
output regulation problems without Nussbaum-type gains
\cite{Nussbaum1983,WangGuay2023,WangGuayBraatz2024,WangGuayXuDochain2025} and
to optimize mechanical systems that are open-loop unstable
\cite{Suttner2023,Suttner2024}. Robustness of large-amplitude, high-frequency
schemes to bounded input and output disturbances is established in
\cite{SuttnerDashkovskiy2022,Labar2022}, and further extensions address
time-varying costs \cite{Grushkovskaya2017}, vanishing input oscillations
\cite{AbdelgalilTaha2021}, unbiased convergence \cite{Yilmaz2024}, and
initialization-free designs \cite{Abdelgalil2024}. In all of these results the
plant is deterministic: disturbances, when present, are bounded functions of
time, and the averaged system is an ordinary differential equation.

Stochastic extremum seeking has developed along two lines. In the first, the
probing signal is random. Manzie and Krsti\'c \cite{ManzieKrstic2009} replace
sinusoidal dithers by stochastic perturbations, and Liu and Krsti\'c
\cite{LiuKrstic2010,LiuKrsticInfinite2010,LiuKrstic2012,LiuKrstic2016} develop
stochastic averaging theorems for fast ergodic excitation and use them to
analyze extremum seeking for static maps and nonlinear dynamical systems, in
continuous and discrete time. In the second, the cost measurement is corrupted
by noise, and convergence is obtained with vanishing gains and stochastic
approximation arguments \cite{Stankovic2010,Khong2015}, through
quasi-stochastic approximation \cite{LauandMeyn2025}, or through time-delay
approaches to averaging that yield mean-square bounds \cite{Yang2026,Zhao2026}.
In both lines the randomness enters through the algorithm or the sensor, the
plant itself is deterministic, and the averaged system is again an ordinary
differential equation. Persistent process noise acting on the plant state is a
different situation. The noise cannot be averaged out, the closed loop is a
diffusion whose stability must be assessed with stochastic Lyapunov tools
\cite{Kushner1967,Khasminskii2012,FreidlinWentzell2012}, and any averaging of
the probing signal must retain the Brownian term
\cite{Khasminskii1968,PavliotisStuart2008}.

Two recent contributions bring Brownian motion into the Lie-bracket framework.
Scheinker \cite{Scheinker2026Ito} extends bounded extremum seeking to It\^o
processes, A cosine-driven stochastic differential equation is compared with a
target diffusion driven by the same Brownian motion on compact time intervals,
with application to gradient-free sampling from generative models. Abdelgalil,
Krsti\'c, and Poveda \cite{Abdelgalil2025,Abdelgalil2026} introduce stochastic
Lie-bracket approximations in which smooth interpolations of a Wiener process
replace periodic dithers. The resulting flows approximate Riemannian Langevin
dynamics on compact manifolds in mean square, and the Gibbs invariant measure
yields infinite-horizon probabilistic guarantees for global optimization. In
both formulations the Brownian term is an ingredient of the optimization
algorithm rather than a disturbance of the plant. In \cite{Scheinker2026Ito}
probing and Brownian forcing share the same coefficient functions, and in
\cite{Abdelgalil2026} the interpolated Wiener process is itself the dither. The
objective is zeroth-order optimization or sampling over the whole state space,
and neither stabilization of an open-loop unstable plant nor compensation of an
unknown constant input arises.

The present paper considers unknown input-affine plants subject to persistent
Brownian disturbances whose diffusion coefficient is independent of the control
input and need not share the actuation directions. The controller uses a single
scalar cost measurement. The principal contribution is a dual-mode, PI
formulation of extremum seeking for this class of controlled diffusions. A
bounded cosine feedback generates, through Lie-bracket averaging, a
proportional action that stabilizes a feasible minimum of the cost and an
integral action that regulates the input to an unknown bias. The averaged closed loop is a PI-regulated diffusion that retains
the process noise, and the analysis is carried out directly on the stochastic
closed loop at a fixed, sufficiently large excitation frequency. Specifically:
\begin{enumerate}
\item[(i)] The averaged dynamics of the cosine-controlled diffusion are
identified as a PI target diffusion. A local matrix inequality condition for
exponential stability of its equilibrium is derived in which the input map
enters only through a positive semidefinite product, so that the proportional
gain can be selected without knowledge of the control direction.
\item[(ii)] For the noise-free closed loop, an exact periodic response is
constructed near the optimum and a Lyapunov function centered on this response
is obtained with bounds that are uniform in the excitation frequency. This
construction separates noise-induced deviations from the oscillations caused by
probing. The mean-square regulation error is bounded by an exponentially
decaying transient plus a term of order $\sigma^{2}+\omega^{-1}$, where
$\sigma$ scales the noise and $\omega$ is the excitation frequency.
\item[(iii)] Capture and residence estimates follow from stopped Lyapunov
arguments \cite{Kushner1967,Khasminskii2012}. For sufficiently small noise,
trajectories starting in a certified compact set enter a prescribed
neighborhood of the optimum and remain there over horizons that grow
exponentially with the inverse noise variance, with a failure probability that
decreases on the same scale, in the spirit of the exit problem of
\cite{FreidlinWentzell2012}.
\item[(iv)] Under an additional regional incremental stability condition of
stochastic contraction type \cite{Pham2009}, the closed loop and the averaged
diffusion driven by the same Brownian motion are compared uniformly in time,
with a mean-square error of order $\omega^{-1}$ that applies until either
trajectory leaves the comparison region; see \cite{Crisan2026} for
uniform-in-time averaging of diffusions under strong exponential stability.
\end{enumerate}
A planar hydrodynamic saddle trap illustrates the regulation result. The
controlled extensional direction is unstable, the passive compressional
direction carries Brownian noise, and an unknown offset requires a nonzero
holding input. A regional target certificate covers uncertain strain rate,
optical amplitude, drift offset, and either actuator sign. The scalar measured
objective is a Gaussian optical-intensity score.

Sections~\ref{secII} and~\ref{secIII} introduce the preliminaries, plant model,
and stability assumptions. Section~\ref{secIV} presents the cosine controller,
and Section~\ref{secV} derives the averaged proportional--integral dynamics and
examines their stability. Section~\ref{secVI} constructs the periodic change of
coordinates, establishes the finite-horizon approximation, and proves
attraction to the noise-free periodic response. Section~\ref{secVII} develops
the stochastic regulation, capture, and residence estimates, followed by the
time-uniform approximation result under incremental stability.
Section~\ref{secVIII} presents the hydrodynamic-trap model and stochastic
simulations, and Section~\ref{secIX} concludes the paper.
Appendices~\ref{app:trap-certificate} and~\ref{app:trap-passive} provide
the regional certificate and an analytical bound on radial residence.

\section{Preliminaries}\label{secII}
The symbols $|\cdot|$, $\|\cdot\|$, and $\|\cdot\|_F$ denote Euclidean, induced, and Frobenius norms. Write $\sym A=(A+A\transpose)/2$, let $\Spp^m$ denote the symmetric positive-definite matrices, and set $\Torus=\R/(2\pi\mathbb Z)$. For open sets, $U\Subset D$ means that $\overline U$ is compact and contained in $D$. A $C^j$ norm includes spatial derivatives through order $j$. Variables entering trigonometric functions are dimensionless. Constants denoted by $C$ may be enlarged between estimates. Their dependence on a time horizon or a noise bound is indicated when relevant.

All stochastic equations are interpreted in the It\^o sense on a filtered probability space satisfying the usual conditions. The Brownian dimension is denoted by $d_B$, reserving $r$ for a Lyapunov level. Solutions are considered up to their maximal lifetime. For an open set $U$, define $\tau_U(Z)=\inf\{t\ge0:Z_t\notin U\}$, with $\inf\varnothing=\infty$. An integrand multiplied by $\ind_{\{t<\tau\}}$ is defined to be zero after exit. For vector fields $A$, $B$ and scalar function $H$, the Lie derivative and the Lie bracket are denoted by
\begin{equation}\label{eq:bracket-convention}
 L_AH=\nabla H\transpose A,\qquad [A,B]=DB\,A-DA\,B
\end{equation}
respectively. 

For a vector- or matrix-valued function $F(\theta,z)$
that is periodic in $\theta$, we use the spatial norm
\[
 \|F\|_{C^k}
 :=
 \max_{0\le j\le k}
 \sup_{\substack{\theta\in\Torus\\z\in\R^{n+m}}}
 \|D_z^jF(\theta,z)\|.
\]
Here $D_z^0F=F$. For $j=0$, the norm is the
Euclidean norm for vectors and the Frobenius norm
for matrices.  For $j\ge1$, it is the corresponding
induced multilinear norm. Thus the $C^k$-norm
bounds the function and its spatial derivatives
through order $k$, uniformly in the phase.
Derivatives with respect to $\theta$ are not
included.
\subsection{A stopped Lyapunov estimate}
Consider the It\^o diffusion $\dd Z_t=a(t,Z_t)\dd t+\sigma S(t,Z_t)\dd B_t$ on an open $U\subset\R^d$, from a deterministic state $z_0$. The coefficients are continuous and locally Lipschitz in the state, uniformly on finite time intervals. For a nonnegative $\Phi\in C^{1,2}([0,\infty)\times U)$, the time-augmented generator $ \cL_\sigma$ acts as
\[
 \cL_\sigma\Phi=\partial_t\Phi+\nabla_z\Phi\transpose a
       +\tfrac12\sigma^2\tr(SS\transpose D_z^2\Phi).
\]
Setting $Q_t:=\Phi(t,Z_t)$, It\^o's formula gives
\[
dQ_t = (\cL_\sigma \Phi)(t, Z_t)\dd t + \sigma \nabla_z\Phi(t, Z_t)\transpose S(t, Z_t)\dd B_t.
\]
Here, $q_\sigma(\Phi):=\sigma^2|S\transpose\nabla_z\Phi|^2$ denotes the quadratic-variation rate of the local-martingale part of $Q_t$. Suppose that the closed sublevels $\{z\in U:\Phi(t,z)\le R\}$ lie in one compact subset of $U$, uniformly in $t$, and let $\tau_R=\inf\{t\ge0:Q_t\ge R\}$. Then for $Q_0 < R$, compact continuation and continuity give $Q_{\tau_R}=R$ at every finite exit.

The following lemma provides a moment bound on $Q_t$ before exit from a Lyapunov sublevel and a bound on the probability of exit.

\begin{lemma}\label{lem:stopped}
Suppose, for $\Phi<R$,
\[
\cL_\sigma\Phi\le-\lambda\Phi+c_0\sigma^2,
 \qquad q_\sigma(\Phi)\le d_0\sigma^2\Phi,
\]
where $\lambda>0$ and $c_0,d_0\ge0$. For $Q_0<R$,
\begin{equation}\label{eq:abstract-moment}
 \E[Q_t\ind_{\{t<\tau_R\}}]
 \le e^{-\lambda t}Q_0+\frac{c_0\sigma^2}{\lambda}(1-e^{-\lambda t}).
\end{equation}
There are $\eta,C_0>0$, depending only on $\lambda,c_0,d_0$, such that for $\sigma>0$ and finite $T\ge0$,
\begin{equation}\label{eq:abstract-exit}
 \Prob\{\tau_R\le T\}
 \le e^{-\eta(R-Q_0)/\sigma^2}+C_0T e^{-\eta R/\sigma^2}.
\end{equation}
For $\sigma=0$, $\tau_R=\infty$ and $Q_t\le e^{-\lambda t}Q_0$.
\end{lemma}
\begin{proof}
It\^o's formula applied up to $\tau_R$ for $e^{\lambda t}Q_t$ yields
\[
\begin{split}
 &e^{\lambda t}\E[Q_t\ind_{\{t<\tau_R\}}]
 +R\E[e^{\lambda\tau_R}\ind_{\{\tau_R\le t\}}]\\
 &\qquad\le Q_0+\frac{c_0\sigma^2}{\lambda}(e^{\lambda t}-1).
\end{split}
\]
Dropping the exit contribution proves \eqref{eq:abstract-moment}. For \eqref{eq:abstract-exit}, choose $\eta>0$ with $\eta d_0\le\lambda$ and set $F=\exp(\eta\Phi/\sigma^2)$. The generator chain rule gives
\[
 \cL_\sigma F\le e^{\eta y}(\eta c_0-\eta\lambda y/2),
 \qquad y=\Phi/\sigma^2.
\]
The right side has a finite upper bound $C_0>0$ on $y\ge0$, independent of $R,\sigma$. Thus $\E F(T\wedge\tau_R,Z_{T\wedge\tau_R})\le e^{\eta Q_0/\sigma^2}+C_0T$. On exit $F=e^{\eta R/\sigma^2}$, which proves \eqref{eq:abstract-exit}. All stopped stochastic integrals are integrable on finite intervals by compactness. For $\sigma=0$, the differential inequality prevents $Q$ from reaching $R$, and compact continuation gives the claim.
\end{proof}

\section{Problem Formulation}\label{secIII}
Consider the system
\[
\begin{aligned}
 \dd X_t&=[f(X_t)+G(X_t)u_t]\dd t+\sigma\Sigma_0(X_t)\dd B_t,\\
 Y_t&=h(X_t),
\end{aligned}
\]
where $X_t\in D\subset\R^n$, $u_t\in\R^m$, and $Y_t\in\R$. The unknown maps are $f:D\to\R^n$, $G:D\to\R^{n\times m}$, $\Sigma_0:D\to\R^{n\times d_B}$, and $h:D\to\R$. The scalar $\sigma\ge0$ scales the diffusion amplitude and is fixed for a given plant. The controller uses only $Y$ and its own states. %

\begin{assumption}\label{ass:smooth}
The domain $D$ is open, $1\le m\le n$, and $f,G,\Sigma_0,h$ are smooth on $D$.
\end{assumption}
\begin{assumption}\label{ass:min}
There exists $(x^\star,u^\star)\in D\times\R^m$ such that
\[
f(x^\star)+G(x^\star)u^\star=0.
\]
The cost satisfies $\nabla h(x^\star)=0$ and $\nabla^2h(x^\star)\succ0$, and $G(x^\star)$ has full column rank.
\end{assumption}

The objective is the regulation of the optimal cost $h^\star=h(x^\star)$ at the corresponding equilibrium pair $(x^\ast,\,u^\ast)$, despite persistent process noise. These disturbances continually perturb the state, requiring ongoing feedback correction to maintain operation near the desired equilibrium. We therefore seek to derive quantitative bounds on the regulation error and on the probability of reaching and remaining in a prescribed neighborhood of the desired operating point.

We define the input sensitivity
\[
 \ell(x):=G(x)\transpose\nabla h(x)\in\R^m.
\]
Its components are the directional derivatives of the cost along the columns of $G$. For gains $K_P,K_I\in\Spp^m$, we introduce the bias $u_b\in\R^m$, the augmented state $z=\col(x,u_b)$, the desired operating point $z^\star=\col(x^\star,u^\star)$, and the domain $\cD=D\times\R^m$. The deterministic PI target field is $\dot z=F_T(z)$, where
\begin{equation}\label{eq:target}
 F_T(z)=\begin{bmatrix}
 f(x)+G(x)(u_b-K_P\ell(x))\\-K_I\ell(x)
 \end{bmatrix}.
\end{equation}
Here, $F_T$ specifies the desired averaged dynamics. The controller in Section \ref{secIV} uses the measured cost to produce PI actions through Lie-bracket averaging, without evaluating $\ell(x)$.

\begin{definition}\label{def:certificate}
A regional target certificate for a compact initial set $\mathcal K\subset\cD$ consists of a bounded open $\mathcal O\Subset\cD$ containing $z^\star$, a function $V\in C^4$ on a neighborhood of $\overline{\mathcal O}$, constants $c_1,c_2,c_3>0$, and levels $0<r_0<R_0$. For $z\in\mathcal O$,
\begin{equation}\label{eq:certificate-positive}
 c_1|z-z^\star|^2\le V(z)\le c_2|z-z^\star|^2,
\end{equation}
\begin{equation}\label{eq:certificate}
 \nabla V(z)\transpose F_T(z)\le-c_3|z-z^\star|^2.
\end{equation}
The set $\mathcal K$ lies in $\{z\in\mathcal O:V(z)\le r_0\}$, and $\mathcal S_{R_0}:=\{z\in\mathcal O:V(z)\le R_0\}$ is compact and contained in $\mathcal O$.
\end{definition}

Definition~\ref{def:certificate} concerns the deterministic
averaged system $\dot z=F_T(z)$. This certificate is used
below to construct the periodic Lyapunov function
$\Phi_\omega$ for the cosine-controlled system
\eqref{eq:complete}, where $\omega>0$ denotes the controller's base angular excitation frequency. The stochastic drift and
quadratic-variation estimates required by
Lemma~\ref{lem:stopped} are then derived from this
construction.

\begin{assumption}\label{ass:certificate}
Let $\cB\subset\cD$ be an open neighborhood of $z^\star$. For every compact $\mathcal K\subset\cB$, gains $K_P,K_I\in\Spp^m$ can be chosen so that $F_T$ admits a regional target certificate for $\mathcal K$. The gains may depend on $\mathcal K$ but are independent of $\omega$.
\end{assumption}

\begin{remark}\label{rem:target-certificate}
Assumption~\ref{ass:certificate} is motivated by earlier work on  deterministic
stability analysis of PI and dual-mode
extremum seeking \cite{Guay2016,GuayAtta2018}.
For certain classes of nonlinear systems with convex cost functions,
these works develop stabilizing gain choices and Lyapunov
constructions under suitable conditions on the dynamics
and input sensitivity. Their analyses address both convergence of the plant state to the optimum and
compensation of the unknown sustaining input.

For the considered plant class, the present formulation expresses the required deterministic
stability property through the regional certificate in
Definition~\ref{def:certificate}. Proposition~\ref{prop:primitive}
provides a local sufficient condition for the existence of this certificate. The subsequent analysis uses this certificate to quantify the effects of periodic excitation and process noise on the corresponding closed loop. 
\end{remark}

\section{Dual-Mode Controller} \label{secIV}
We choose $\omega>0$, distinct positive integers $\kappa_1,\ldots,\kappa_m$, and set $C_P=K_P^{1/2}$ and $C_I=K_IK_P^{-1/2}$. For $v_t^\omega=\col(v_{1,t}^\omega,\ldots,v_{m,t}^\omega)$, the feedback is
\begin{equation}\label{eq:controller}
\begin{aligned}
 \dd\theta_t&=\omega\dd t,\qquad\theta_0\in\Torus,\\
 v_{i,t}^\omega&=\sqrt{2\kappa_i\omega}\cos(Y_t^\omega+\kappa_i\theta_t),\\
 u_t^\omega&=u_{b,t}^\omega+C_Pv_t^\omega,\\
 \dd u_{b,t}^\omega&=C_Iv_t^\omega\dd t.
\end{aligned}
\end{equation}
Here $Y_t^\omega=h(X_t^\omega)$ and $u_{b,0}^\omega$ is prescribed. The phase and bias have finite-variation paths. The bias is random through $Y_t^\omega$, but its update requires neither an output derivative nor integration against $\dd Y_t^\omega$.

For $z=\col(x,u_b)$, we define the augmented drift
\[
 b(z)=\begin{bmatrix}f(x)+G(x)u_b\\0_m\end{bmatrix}
\]
and the diffusion matrix
\[
 S_0(z)=\begin{bmatrix}\Sigma_0(x)\\0_{m\times d_B}\end{bmatrix}.
\]
With $e_i$ the $i$th standard basis vector of $\R^m$, we introduce the following augmented-state fields
\begin{equation}\label{eq:fields}
\begin{aligned}
 A_i(z)=&{}\cos(h(x))H_i(z)\\
 B_i(z)=&-\sin(h(x))H_i(z),
\end{aligned}
\end{equation}
where \[ H_i(z)=\begin{bmatrix}G(x)C_Pe_i\\C_Ie_i\end{bmatrix}.\]
 Let $p_i(\theta)=\sqrt{2\kappa_i}\cos(\kappa_i\theta)$ and $q_i(\theta)=\sqrt{2\kappa_i}\sin(\kappa_i\theta)$, and put $a(\theta,z)=\sum_i[A_i(z)p_i(\theta)+B_i(z)q_i(\theta)]$. The augmented closed loop is
\begin{equation}\label{eq:complete}
 \dd Z_t^\omega=[b(Z_t^\omega)+\sqrt\omega\,a(\theta_t,Z_t^\omega)]\dd t
       +\sigma S_0(Z_t^\omega)\dd B_t.
\end{equation}
Its phase-augmented generator $\cA_{\omega,\sigma}$ is defined such that
\begin{equation}\label{eq:extended-generator}
\begin{aligned}
 \cA_{\omega,\sigma}\Phi
 ={}\omega\partial_\theta\Phi+&\nabla_z\Phi\transpose(b+\sqrt\omega a)\\
       &+\tfrac12\sigma^2\tr(S_0S_0\transpose D_z^2\Phi).
\end{aligned}
\end{equation}

With $a_\omega=(2\omega\sum_i\kappa_i)^{1/2}$, the bounded cosine yields the following bounds
\[
|u_t^\omega-u_{b,t}^\omega|\le\|C_P\|a_\omega,
 \qquad |C_Iv_t^\omega|\le\|C_I\|a_\omega.
\]
For each fixed frequency, the bounds are independent of the cost magnitude and grow with $\sqrt\omega$. Replacing the amplitude in \eqref{eq:controller} by a frequency-independent constant $\alpha_i$ results in the bracket coefficient $\alpha_i^2/(2\kappa_i\omega)$. Hence, the chosen amplitude scaling $\sqrt\omega$ maintains a non-vanishing Lie-bracket coefficient as the frequency increases.

\begin{remark}\label{rem:washout}
A constant cost offset leaves $\ell(x)$ and $F_T$ unchanged but generally changes the feedback \eqref{eq:controller}. To make the applied input invariant to such offsets, replace $Y_t^\omega$ in \eqref{eq:controller} by $Y_t^\omega-w_t$ with $\dd w_t=\gamma(Y_t^\omega-w_t)\dd t$, initialized $w_0=Y_0^\omega$, and fixed $\gamma>0$. This modification adds filter dynamics while preserving the averaged dynamics $F_T$ of $(x, u_b)$. Applying the regional stability results to this modified controller requires a Lyapunov certificate that also includes the filter state $w$.
\end{remark}

\section{Averaged Dynamics and Target Stability}\label{secV}
We first identify the averaged PI dynamics and give a local sufficient condition for the target certificate. For a $2\pi$-periodic function, write $\avg{v}=(2\pi)^{-1}\int_0^{2\pi}v(\theta)\dd\theta$. The normalized input dithers satisfy
\[
 \left\langle\left(\int_0^\theta p_i(s)\dd s\right)q_i(\theta)\right\rangle=1.
\]
Reversing $p_i,q_i$ gives $-1$, while all corresponding averages between distinct channels vanish.

\begin{proposition}\label{prop:bracket}
The fields in \eqref{eq:fields} satisfy
\begin{equation}\label{eq:PI-bracket}
 \sum_{i=1}^m[A_i,B_i]
 =-\begin{bmatrix}GK_P\ell\\K_I\ell\end{bmatrix}.
\end{equation}
\end{proposition}
\begin{proof}
For $\widetilde h(z)=h(x)$, the product rule with \eqref{eq:bracket-convention} gives $[A_i,B_i]=-(L_{H_i}\widetilde h)H_i$. Since $L_{H_i}\widetilde h=e_i\transpose C_P\transpose\ell$, summation gives
\[
 \sum_i[A_i,B_i]=-
 \begin{bmatrix}GC_PC_P\transpose\ell\\ C_IC_P\transpose\ell\end{bmatrix}.
\] The identities $C_PC_P\transpose=K_P$ and $C_IC_P\transpose=K_I$ prove \eqref{eq:PI-bracket}.
\end{proof}
We note that adding the bracket correction in \eqref{eq:PI-bracket} to the drift $b$ gives $F_T$. The associated averaged diffusion is
\begin{equation}\label{eq:average-sde}
 \dd\bar Z_t=F_T(\bar Z_t)\dd t+\sigma S_0(\bar Z_t)\dd B_t.
\end{equation}

Setting $\sigma=0$ in \eqref{eq:average-sde} recovers the deterministic PI target $\dot z = F_T(z)$. By Definition \ref{def:certificate}, if a regional target certificate is available, then \eqref{eq:certificate-positive}, \eqref{eq:certificate} imply that $\mathcal S_{R_0}$ is forward invariant, and every solution initialized in $\mathcal S_{R_0}$ satisfies
\[
 |z(t)-z^\star|\le\sqrt{c_2/c_1}\,e^{-c_3t/(2c_2)}|z(0)-z^\star|
\]
Compact continuation gives existence for all $t\ge0$. The following proposition gives a sufficient condition for such a certificate to exist locally.

\begin{proposition}\label{prop:primitive}
Under Assumptions~\ref{ass:smooth}--\ref{ass:min}, define $H_\star=\nabla^2h(x^\star)$, $G_\star=G(x^\star)$, and $A_\star=D_x(f+Gu^\star)(x^\star)$. If $K_P\in\Spp^m$ satisfies
\begin{equation}\label{eq:primitive}
 \sym(H_\star A_\star)-H_\star G_\star K_PG_\star\transpose H_\star\prec0,
\end{equation}
then $z^\star$ is locally exponentially stable for every fixed $K_I\in\Spp^m$. Every compact initial set in a sufficiently small neighborhood admits a regional target certificate.
\end{proposition}
\begin{proof}
Set $A_P=A_\star-G_\star K_PG_\star\transpose H_\star$. Since $D\ell(x^\star)=G_\star\transpose H_\star$, the target Jacobian $A_T:=DF_T(z^\star)$ is then
\[
 A_T=\begin{bmatrix}A_P&G_\star\\-K_IG_\star\transpose H_\star&0\end{bmatrix}.
\]
By \eqref{eq:primitive}, $Q=-\sym(H_\star A_P)\succ0$. Then defining $P_0=\diag(H_\star,K_I^{-1})$, we have $A_T\transpose P_0+P_0A_T=\diag(-2Q,0)$. An eigenpair $A_Tv=\mu v$, $v=\col(v_x,v_b)\ne0$, satisfies
\[
 2\operatorname{Re}(\mu)v^*P_0v=-2v_x^*Qv_x.
\]
If $\operatorname{Re}(\mu)\ge0$, then $v_x=0$ and the first block gives $G_\star v_b=0$. Full column rank contradicts $v\ne0$, so $A_T$ is Hurwitz. Consequently, there exists a unique  $P\in\Spp^{m+n}$ such that $A_T\transpose P+PA_T=-I$. Denoting \(e:=z-z^\ast\), smoothness gives
\[
 \frac{\dd}{\dd t}(e\transpose Pe)\le-\tfrac12|e|^2
\]
on a sufficiently small closed ball. We choose an outer quadratic sublevel strictly inside that ball. Its interior contains a neighborhood of $z^\star$, and each compact subset lies below a smaller positive level. These choices give Definition~\ref{def:certificate}, with $c_1=\lambda_{\min}(P)$, $c_2=\lambda_{\max}(P)$, and $c_3=1/2$.
\end{proof}
For fixed $(H_\star, G_\star, A_\star)$, conditions \eqref{eq:primitive} and $K_P\succ 0$ form a strict LMI feasibility problem in $K_P$, providing a sufficient condition for local exponential stability of the deterministic PI target.

Fix a regional target certificate from
Definition~\ref{def:certificate}, with outer sublevel
$\mathcal S_{R_0}=\{z\in\mathcal O:V(z)\le R_0\}$.
For each $r_0<R\le R_0$, define
\[
 \mathcal S_R:=\{z\in\mathcal O:V(z)\le R\}.
\]
This is a closed subset of the compact set
$\mathcal S_{R_0}$ and is therefore compactly
contained in $\mathcal O$. Moreover,
$\mathcal K\subset\{z\in\mathcal O:V(z)<R\}$.

For the averaged diffusion \eqref{eq:average-sde},
the time-augmented generator applied to $V$ is
\[
 (\overline{\cL}_\sigma V)(z)
 =
 \nabla V(z)\transpose F_T(z)
 +\frac{\sigma^2}{2}
 \tr\!\left(S_0(z)S_0(z)\transpose D^2V(z)\right).
\]
The quadratic bounds on $V$ from \eqref{eq:certificate-positive} imply
$\nabla V(z^\star)=0$. A local Hessian bound near
$z^\star$, together with compactness away from
$z^\star$, gives a constant $L_V$ such that
\[
 |\nabla V(z)|\le L_V|z-z^\star|
\]
for \(z\in\mathcal S_{R_0}.\)
The target certificate and the definition of
$q_\sigma(V)$ consequently give
\begin{equation}\label{eq:average-generator-bound}
\begin{aligned}
 \overline{\cL}_\sigma V
 &\le-\frac{c_3}{c_2}V+C_V\sigma^2,\\
 q_\sigma(V)
 &\le D_V\sigma^2V,
\end{aligned}
\end{equation}
throughout $\mathcal S_{R_0}$. One may take
\[
 C_V
 :=
 \frac12\max_{z\in\mathcal S_{R_0}}
 \left(\|S_0(z)\|_F^2\|D^2V(z)\|\right)
\]
and
\[
 D_V
 :=
 \frac{L_V^2}{c_1}
 \max_{z\in\mathcal S_{R_0}}\|S_0(z)\|_F^2.
\]
Both constants are finite and independent of
$\sigma$ and the chosen level $R$.

For $\bar Z_0=z_0\in\mathcal K$, let
\[
 \bar\tau_R
 :=
 \inf\{t\ge0:V(\bar Z_t)\ge R\}
\]
be the first exit from $\{z\in\mathcal O:V(z)<R\}$.
Lemma~\ref{lem:stopped} applies to
$Q_t=V(\bar Z_t)$ up to $\bar\tau_R$, using
\eqref{eq:average-generator-bound}. To obtain corresponding estimates for the augmented closed loop \eqref{eq:complete}, we must also account for the oscillations generated by the probing signal.
\section{Periodic Normal Form and Deterministic Closed-Loop Stability}
\label{secVI}
A periodic change of coordinates removes the leading oscillatory drift. We use the resulting perturbation of $F_T$ to construct an exact periodic response, then center the target Lyapunov function on that response.

\subsection{A periodic change of coordinates}
\label{sec:normal-coordinate}
Fix a certificate and choose bounded open sets with
\[
 \mathcal S_{R_0}\subset\mathcal O_1\Subset\mathcal O_2\Subset\mathcal O.
\]
Choose a smooth cutoff function
$\chi:\R^{n+m}\to[0,1]$ that equals one on a
neighborhood of $\overline{\mathcal O_2}$ and has
compact support in $\mathcal O$.
For each $F\in\{b,S_0,A_i,B_i\}$, define
\[
 \widetilde F(z)
 :=
 \begin{cases}
  \chi(z)F(z), & z\in\mathcal O,\\
  0, & z\notin\mathcal O.
 \end{cases}
\]

Define the averaged drift of the extended system by
\[
 \widetilde F_T(z)
 :=
 \widetilde b(z)
 +\sum_{i=1}^m
 [\widetilde A_i,\widetilde B_i](z).
\]
Since the cutoff equals one near
$\overline{\mathcal O_2}$, the extended fields and
their derivatives agree with the original ones there.
In particular,
\[
 \widetilde F_T(z)=F_T(z),
 \qquad z\in\mathcal O_2,
\]
where $F_T$ is the deterministic target drift
defined in \eqref{eq:target}.

We construct the coordinate transformation using
the extended coefficients, whose spatial derivatives
are globally bounded. The resulting identities apply
to the original system whenever the physical and
transformed states remain in $\mathcal O_2$.

For the normal-form construction in this subsection
and the auxiliary comparison in the proof of
Proposition~\ref{prop:finite-window}, we suppress
the tildes and write $b,S_0,A_i,B_i$ for the
compactly supported extensions. Within these
calculations, $a=\sum_i(A_ip_i+B_iq_i)$ and
$F_T=b+\sum_i[A_i,B_i]$ are formed from the
extended fields.

Outside these calculations, the symbols retain
their original definitions. The two versions
agree on $\mathcal O_2$, where the cutoff equals
one. Thus the identities obtained for the extended
equations apply to the original equations whenever
the relevant states remain in $\mathcal O_2$.

Throughout, we set $\eps=\omega^{-1/2}$. Then \eqref{eq:complete} has drift $b+\eps^{-1}a$ and $\dd\theta_t=\eps^{-2}\dd t$.  We denote the zero-mean periodic corrector as
\[
 U_1(\theta,z)=\sum_i\sqrt{2/\kappa_i}
 [A_i(z)\sin(\kappa_i\theta)-B_i(z)\cos(\kappa_i\theta)].
\]
It satisfies $\partial_\theta U_1=a$. Let $U_2$ be the unique zero-mean periodic solution of
\begin{equation}\label{eq:homological}
 \partial_\theta U_2=D_za\,U_1-\avg{D_za\,U_1},
\end{equation}
and define
\[
\Gamma_\eps(\theta,\zeta)=\zeta+\eps U_1(\theta,\zeta)+\eps^2U_2(\theta,\zeta).
\]
All norms in the next lemma are global spatial norms for the extensions, uniform in phase.

\begin{lemma}\label{lem:normal}
There exist $\eps_0,C>0$ such that $\Gamma_\eps(\theta,\cdot)$ is a diffeomorphism for $0<\eps\le\eps_0$ and every phase, with
\begin{equation}\label{eq:gamma-bounds}
 \|\Gamma_\eps-\mathrm{id}\|_{C^3}
 +\|\Gamma_\eps^{-1}-\mathrm{id}\|_{C^3}\le C\eps.
\end{equation}
For the extended closed loop, $\zeta_t=\Gamma_\eps^{-1}(\theta_t,Z_t^\omega)$ satisfies
\begin{equation}\label{eq:normal-form}
\begin{split}
 \dd\zeta_t={}&[F_T(\zeta_t)+\eps R_\eps(\theta_t,\zeta_t)
                   +\sigma^2B_\eps(\theta_t,\zeta_t)]\dd t\\
 &+\sigma\widehat S_\eps(\theta_t,\zeta_t)\dd B_t.
\end{split}
\end{equation}
The remainder fields are periodic in phase, independent of $\sigma$, and satisfy $\|R_\eps\|_{C^2}\le C$ and
\begin{equation}\label{eq:normal-bounds}
 \|B_\eps\|_{C^2}+\|\widehat S_\eps-S_0\|_{C^2}\le C\eps.
\end{equation}
The constants are independent of $\eps,\theta,\sigma$. For the original coefficients, \eqref{eq:normal-form} holds until $Z_t^\omega$ exits $\mathcal O_1$.
\end{lemma}
\begin{proof}
The probing-signal identities give $\avg{D_za\,U_1}=\sum_i[A_i,B_i]=F_T-b$. The smooth extensions give bounded spatial derivatives of $U_1$ through order six and $U_2$ through order five. Put $P_\eps=U_1+\eps U_2$ and $J_\eps=D_\zeta\Gamma_\eps=I+\eps D_\zeta P_\eps$. Choose $\eps_0$ so that $\eps\|D_\zeta P_\eps\|_\infty\le1/2$. The equation $\zeta=z-\eps P_\eps(\theta,\zeta)$ is a contraction on $\R^{n+m}$, defining the inverse. Differentiating this identity proves \eqref{eq:gamma-bounds} and gives
\begin{equation}\label{eq:bilipschitz}
 \tfrac12|v-w|\le|\Gamma_\eps(\theta,v)-\Gamma_\eps(\theta,w)|
 \le\tfrac32|v-w|.
\end{equation}

With unmarked arguments evaluated at $(\theta,\zeta)$, Taylor's formula and \eqref{eq:homological} give
\[
 b(\Gamma_\eps)+\eps^{-1}a(\theta,\Gamma_\eps)
 -\eps^{-2}\partial_\theta\Gamma_\eps=F_T+\eps Q_\eps,
\]
where
\[
\begin{split}
 Q_\eps={}&\int_0^1Db(\zeta+s\eps P_\eps)P_\eps\dd s+D_\zeta a\,U_2\\
 &+\int_0^1(1-s)D_\zeta^2a(\theta,\zeta+s\eps P_\eps)
                   [P_\eps,P_\eps]\dd s.
\end{split}
\]
The two corrector equations cancel the $\eps^{-1}$ and order-one oscillatory terms. The integral formula gives a uniform $C^2$ bound for $Q_\eps$. Setting
\[
 R_\eps=J_\eps^{-1}[Q_\eps-D_\zeta P_\eps F_T]
\]
therefore gives $J_\eps^{-1}(F_T+\eps Q_\eps)=F_T+\eps R_\eps$ and $\|R_\eps\|_{C^2}\le C$.

It\^o's formula for $Z=\Gamma_\eps(\theta,\zeta)$ gives
\[
 \dd Z=\eps^{-2}\partial_\theta\Gamma_\eps\dd t+J_\eps\dd\zeta
 +\frac{\sigma^2}{2}\sum_{j=1}^{d_B}D_\zeta^2\Gamma_\eps
 [\widehat S_{\eps,j},\widehat S_{\eps,j}]\dd t.
\]
Here $\widehat S_{\eps,j}$ is the $j$th column and $D_\zeta^2\Gamma_\eps[v,w]$ is the vector of component-Hessian bilinear forms. Matching diffusion and drift determines
\[
\widehat S_\eps=J_\eps^{-1}S_0(\Gamma_\eps)
\]
and 
\[
 B_\eps=-\tfrac12J_\eps^{-1}\sum_{j=1}^{d_B}
 D_\zeta^2\Gamma_\eps[\widehat S_{\eps,j},\widehat S_{\eps,j}].
\]
Thus $\sigma^2B_\eps$ is the coordinate-induced It\^o correction. The corrector bounds give $J_\eps^{-1}-I=O(\eps)$ and $D_\zeta^2\Gamma_\eps=O(\eps)$ in $C^2$, proving \eqref{eq:normal-bounds}. These formulas contain no $\sigma$.

Finally, the distance from $\overline{\mathcal O_1}$ to $\R^{n+m}\setminus\mathcal O_2$ is positive. We decrease $\eps_0$ until the displacement in \eqref{eq:gamma-bounds} is smaller than this distance. Both $Z_t^\omega$ and $\Gamma_\eps^{-1}(\theta_t,Z_t^\omega)$ then lie in $\mathcal O_2$ before exiting $\mathcal O_1$, where the extensions and original fields agree.
\end{proof}
The averaged diffusion \eqref{eq:average-sde} retains the plant’s process-noise term while replacing the oscillatory feedback by its effective PI drift. The following proposition compares \eqref{eq:average-sde} with the original closed loop \eqref{eq:complete}, up to a common exit time.

\begin{proposition}\label{prop:finite-window}
Under Lemma~\ref{lem:normal}, fix $\mathcal K\subset\mathcal U\Subset\mathcal O_1$ with $\mathcal K$ compact and $\mathcal U$ open. Systems \eqref{eq:complete} and \eqref{eq:average-sde} are subject to the same Brownian motion and common deterministic initial state $z_0$. Set $\tau^\omega=\tau_{\mathcal U}(Z^\omega)\wedge\tau_{\mathcal U}(\bar Z)$. For every finite $T>0$ and $\bar\sigma\ge0$, there exist $\omega_0,C_{\mathcal K,T,\bar\sigma}>0$ such that
\begin{equation}\label{eq:finite-window}
 \sup_{z_0\in\mathcal K,\,\theta_0\in\Torus}
 \E\sup_{t\le T\wedge\tau^\omega}|Z_t^\omega-\bar Z_t|^2
 \le C_{\mathcal K,T,\bar\sigma}/\omega
\end{equation}
for $\omega\ge\omega_0$ and $0\le\sigma\le\bar\sigma$.
\end{proposition}
\begin{proof}
First use the extended versions of \eqref{eq:complete} and \eqref{eq:average-sde}, with solutions $Z^{\omega,\mathrm e},\bar Z^{\mathrm e}$. Put $\zeta_t=\Gamma_\eps^{-1}(\theta_t,Z_t^{\omega,\mathrm e})$ and $e_t=\zeta_t-\bar Z_t^{\mathrm e}$. Subtracting \eqref{eq:average-sde} from \eqref{eq:normal-form} gives $\dd e_t=r_t\dd t+\Xi_t\dd B_t$, where
\[
 r_t=F_T(\zeta_t)-F_T(\bar Z_t^{\mathrm e})
       +\eps R_\eps(\theta_t,\zeta_t)+\sigma^2B_\eps(\theta_t,\zeta_t)
\]
and $\Xi_t=\sigma[\widehat S_\eps(\theta_t,\zeta_t)-S_0(\bar Z_t^{\mathrm e})]$. The extended $F_T,S_0$ are globally Lipschitz, so Lemma~\ref{lem:normal} gives
\[
 |e_0|^2\le C\eps^2,\qquad
 |r_t|^2+\|\Xi_t\|_F^2\le C_{\bar\sigma}(|e_t|^2+\eps^2).
\]
The integral error equation, Cauchy--Schwarz, and Doob's inequality with It\^o's isometry imply, for $\mathcal E(t)=\E\sup_{s\le t}|e_s|^2$ and $t\le T$,
\[
\begin{split}
 \mathcal E(t)&\le3|e_0|^2+3T\int_0^t\E|r_s|^2\dd s
                       +12\int_0^t\E\|\Xi_s\|_F^2\dd s\\
 &\le C_{T,\bar\sigma}\eps^2+C_{T,\bar\sigma}\int_0^t\mathcal E(s)\dd s.
\end{split}
\]
Gr\"onwall gives $\mathcal E(T)\le C_{T,\bar\sigma}\eps^2$. The identity $Z^{\omega,\mathrm e}-\bar Z^{\mathrm e}=\Gamma_\eps(\theta,\zeta)-\zeta+e$ and \eqref{eq:gamma-bounds} give the same order in physical coordinates. The extended and original solutions agree through $\tau^\omega$ by pathwise uniqueness and continuity. Restriction to this interval and $\eps^2=\omega^{-1}$ prove the result.
\end{proof}

\subsection{Periodic response of the noise-free closed loop}
\label{sec:orbit}
We construct a periodic trajectory near $z^\star$ for fixed frequency. It is the reference for the stochastic estimates, rather than an equilibrium of the time-periodic physical system.

\begin{lemma}\label{lem:orbit}
For a certified target and sufficiently small $\eps$, there is a unique $2\pi$-periodic $C^1$ function $\chi_\eps$ in a fixed sufficiently small supremum-norm neighborhood of zero satisfying
\begin{equation}\label{eq:chi}
 \omega\frac{\dd\chi_\eps}{\dd\theta}
 =F_T(z^\star+\chi_\eps)+\eps R_\eps(\theta,z^\star+\chi_\eps).
\end{equation}
It obeys $\|\chi_\eps\|_\infty\le C\eps$. The physical response
\begin{equation}\label{eq:physical-orbit}
 z_\omega^\circ(\theta)=\Gamma_\eps(\theta,z^\star+\chi_\eps(\theta))
\end{equation}
satisfies $\sup_\theta|z_\omega^\circ(\theta)-z^\star|\le C\eps$. For every initial phase, $t\mapsto z_\omega^\circ(\theta_0+\omega t)$ is a $2\pi/\omega$-periodic solution of \eqref{eq:complete} with $\sigma=0$.
\end{lemma}
\begin{proof}
Expanding \eqref{eq:certificate-positive}--\eqref{eq:certificate} at $z^\star$ gives $P=\tfrac12D^2V(z^\star)\succ0$ and $A_T\transpose P+PA_T\preceq-c_3I$, where $A_T=DF_T(z^\star)$. Thus $\|e^{A_Tt}\|\le M_Ae^{-\alpha_A t}$. We define $N(e)=F_T(z^\star+e)-A_Te$. Locally, $|N(e)|\le L_N|e|^2$ and $\|DN(e)\|\le2L_N|e|$. Let $L_R$ bound $R_\eps$ and its first spatial derivative, uniformly in $\eps$.

On continuous periodic functions, we define
\[
\begin{split}
 (\mathcal T_\eps\chi)(\theta)=\int_0^\infty e^{A_Ts}
 \big[&N(\chi(\theta-\omega s))\\
 &+\eps R_\eps(\theta-\omega s,z^\star+\chi(\theta-\omega s))\big]\dd s.
\end{split}
\]
Writing $c_A=M_A/\alpha_A$, we choose $\rho_*>0$ with $\overline B_{2\rho_*}(z^\star)\subset\mathcal O_1$ and $2c_AL_N\rho_*\le1/4$. We can then take $c_A\eps L_R\le\min\{1/4,\rho_*/2\}$. On the closed ball $\|\chi\|_\infty\le\rho_*$, the map has a contraction constant of at most $c_A(2L_N\rho_*+\eps L_R)\le1/2$ and $\|\mathcal T_\eps0\|_\infty\le c_A\eps L_R\le\rho_*/2$. It maps the ball into itself, giving a unique fixed point and $\|\chi_\eps\|_\infty\le2c_AL_R\eps$.

For $q(t)=\chi_\eps(\theta_0+\omega t)$, the fixed-point identity is
\[
 q(t)=\int_{-\infty}^t e^{A_T(t-s)}
 [N(q(s))+\eps R_\eps(\theta_0+\omega s,z^\star+q(s))]\dd s.
\]
The integrand is continuous and bounded after its exponential factor is removed. Differentiation gives \eqref{eq:chi}. Conversely, variation of constants for any bounded complete periodic solution, followed by resetting the starting time to $-\infty$, gives this same identity. This observation proves uniqueness in the fixed ball. Phase translations preserve the supremum norm, so the estimates are uniform in frequency. The displacement bounds keep the response in $\mathcal O_1$, and the exact transformation with $\sigma=0$ proves the assertion for \eqref{eq:complete}. Periodicity follows from that of $\Gamma_\eps,\chi_\eps$.
\end{proof}

\subsection{Centering on the periodic response}
For $\sigma=0$, write $Z_t^{\omega,0}$ for the solution of \eqref{eq:complete} from state $z_0$ and phase $\theta_0$.
\label{sec:centered}
Define
\[
\Psi_\omega(\theta,z)=\Gamma_\eps^{-1}(\theta,z)-\chi_\eps(\theta).
\]
It maps $z_\omega^\circ(\theta)$ to $z^\star$. Along \eqref{eq:complete} with $\sigma=0$, $y_t=\Psi_\omega(\theta_t,Z_t^{\omega,0})$ satisfies
\begin{equation}\label{eq:centered-ode}
 \dd y_t=[F_T(y_t)+E_\eps(\theta_t,y_t)]\dd t,
\end{equation}
where
\[
\begin{aligned}
 E_\eps(\theta,y)={}&F_T(y+\chi_\eps(\theta))-F_T(y)\\
 &-F_T(z^\star+\chi_\eps(\theta))+F_T(z^\star)\\
 &+\eps R_\eps(\theta,y+\chi_\eps(\theta))\\
 &-\eps R_\eps(\theta,z^\star+\chi_\eps(\theta)).
\end{aligned}
\]
The extended $F_T$ has bounded second derivatives and $R_\eps$ has uniformly bounded first derivatives. Therefore $\|D_yE_\eps\|_\infty\le C\eps$. Since $E_\eps(\theta,z^\star)=0$, integration along the segment from $z^\star$ to $y$ gives
\begin{equation}\label{eq:centered-remainder}
 |E_\eps(\theta,y)|\le C\eps|y-z^\star|.
\end{equation}

Set
\begin{equation}\label{eq:Phi}
 \Phi_\omega(\theta,z)=V(\Psi_\omega(\theta,z)).
\end{equation}
The uniform displacement bounds ensure $\Psi_\omega(\theta,\overline{\mathcal O_1})\subset\mathcal O_2$ for small $\eps$, so this function is defined on a neighborhood of $\Torus\times\overline{\mathcal O_1}$. For $0<s\le R_0$, put $\mathcal S_s^{<}=\{y\in\mathcal O:V(y)<s\}$ and define
\begin{equation}\label{eq:tubes}
 \mathcal D_s^\omega(\theta)
 =\{\Gamma_\eps(\theta,y+\chi_\eps(\theta)):y\in\mathcal S_s^{<}\}.
\end{equation}
Their closures lie in a fixed compact subset of $\mathcal O_1$ for small $\eps$. The defining diffeomorphism maps sublevel boundaries to boundaries, where $\Phi_\omega=s$.

\begin{proposition}\label{prop:deterministic-loop}
Fix $r_0<r<R_0$. For sufficiently large $\omega$, with $\lambda=c_3/(2c_2)$,
\begin{equation}\label{eq:det-Phi}
 \cA_{\omega,0}\Phi_\omega\le-\lambda\Phi_\omega
\end{equation}
on $\mathcal D_{R_0}^\omega(\theta)$, and
\begin{equation}\label{eq:orbital-comparability}
 a_1|z-z_\omega^\circ(\theta)|^2\le\Phi_\omega(\theta,z)
 \le a_2|z-z_\omega^\circ(\theta)|^2,
\end{equation}
with $a_1=4c_1/9$ and $a_2=4c_2$. Every $z_0\in\mathcal D_s^\omega(\theta_0)$, $0<s\le R_0$, generates a complete noise-free solution that remains in $\mathcal D_s^\omega(\theta_t)$ and satisfies
\begin{equation}\label{eq:det-orbit-decay}
 |Z_t^{\omega,0}-z_\omega^\circ(\theta_t)|^2
 \le\frac{a_2}{a_1}e^{-\lambda t}|z_0-z_\omega^\circ(\theta_0)|^2.
\end{equation}
Moreover, $\mathcal K\subset\mathcal D_r^\omega(\theta)$ for every phase. The constants and threshold are uniform over the phase and initial state.
\end{proposition}
\begin{proof}
On $\mathcal S_{R_0}$, $|\nabla V(y)|\le L_V|y-z^\star|$. Equations~\eqref{eq:certificate}, \eqref{eq:centered-ode}, and \eqref{eq:centered-remainder} give
\[
 \cA_{\omega,0}\Phi_\omega
 \le-(c_3-CL_V\eps)|y-z^\star|^2\le-\lambda V(y)
\]
for $CL_V\eps\le c_3/2$. Applying \eqref{eq:bilipschitz} to $y+\chi_\eps$ and $z^\star+\chi_\eps$, then using \eqref{eq:certificate-positive}, proves \eqref{eq:orbital-comparability}.

Until first exit from $\mathcal D_s^\omega(\theta_t)$, \eqref{eq:det-Phi} gives $\Phi_\omega(\theta_t,Z_t^{\omega,0})\le e^{-\lambda t}\Phi_\omega(\theta_0,z_0)<s$. At a finite exit, continuity would give $\Phi_\omega=s$, which is a contradiction. The fixed compact containment permits continuation for all time. Applying \eqref{eq:orbital-comparability} at the initial and current states yields \eqref{eq:det-orbit-decay}.

Finally, \eqref{eq:gamma-bounds} and the bound on
$\chi_\eps$ give
\[
 |\Psi_\omega(\theta,z)-z|\le C_\Psi\eps
\]
for \( (\theta,z)\in\Torus\times\mathcal K\)  with $C_\Psi$ independent of $\eps$ and $\theta$.
We choose $\delta>0$ such that the closed neighborhood
$\mathcal K_\delta
:=\{y:\operatorname{dist}(y,\mathcal K)\le\delta\}$
is contained in $\mathcal O$, and let $L_{\mathcal K}$
bound $|\nabla V|$ there. Choose $\eps_1>$ sufficiently small such that
$C_\Psi\eps_1<\delta$ and $L_{\mathcal K}C_\Psi\eps_1<r-r_0$.

For $z\in\mathcal K$, $\theta\in\Torus$, and
$0<\eps\le\eps_1$, the point
$y=\Psi_\omega(\theta,z)$ and the segment joining
it to $z$ lie in $\mathcal K_\delta\subset\mathcal O$.
The gradient bound therefore gives
\[
\begin{aligned}
 \Phi_\omega(\theta,z)=V(y)
 &\le V(z)+L_{\mathcal K}|y-z|\\
 &\le r_0+L_{\mathcal K}C_\Psi\eps<r.
\end{aligned}
\]
Thus $y\in\mathcal S_r^{<}$. Since
$z=\Gamma_\eps(\theta,y+\chi_\eps(\theta))$,
the definition \eqref{eq:tubes} yields
\(\mathcal K\subset\mathcal D_r^\omega(\theta)\)
for every \(\theta\in\Torus\) and \(0<\eps\le\eps_1.
\)
\end{proof}

\section{Stochastic Regulation and Long-Horizon Approximation}
\label{secVII}
The periodic function $\Phi_\omega$ transfers the deterministic decay estimate to the physical diffusion. We first obtain regulation and capture estimates under the target certificate. A separate incremental condition then provides a time-uniform comparison with the bracket diffusion.

\subsection{Regulation about the periodic response}
For \eqref{eq:complete}, define the exit time
\[
 \tau_R^\omega=\inf\{t\ge0:Z_t^\omega\notin\mathcal D_R^\omega(\theta_t)\},
\]
where the sublevels are defined in \eqref{eq:tubes}.

\begin{theorem}\label{thm:main}
Let Assumptions~\ref{ass:smooth}--\ref{ass:certificate} hold. Fix a certificate for $\mathcal K$, levels $r_0<r<R<R_0$, and a finite noise bound $\bar\sigma>0$. There exists $\omega_{\mathcal K}>0$ such that, for every $\omega\ge\omega_{\mathcal K}$, the closed loop \eqref{eq:complete} satisfies
\begin{equation}\label{eq:physical-moment}
\begin{split}
 &\E[|Z_t^\omega-z^\star|^2\ind_{\{t<\tau_R^\omega\}}]\\
 &\le Ce^{-\lambda t}|z_0-z^\star|^2+C(\sigma^2+\omega^{-1})
\end{split}
\end{equation}
where $\lambda=c_3/(2c_2)$. Relative to the response \eqref{eq:physical-orbit},
\begin{equation}\label{eq:orbital-moment}
\begin{split}
 &\E[|Z_t^\omega-z_\omega^\circ(\theta_t)|^2\ind_{\{t<\tau_R^\omega\}}]\\
 &\quad\le Me^{-\lambda t}|z_0-z_\omega^\circ(\theta_0)|^2
                   +C\sigma^2(1-e^{-\lambda t}).
\end{split}
\end{equation}
For $\sigma>0$ and finite $T\ge0$,
\begin{equation}\label{eq:main-exit}
 \Prob\{\tau_R^\omega\le T\}
 \le e^{-\eta(R-r)/\sigma^2}+C_0T e^{-\eta R/\sigma^2}.
\end{equation}
 The positive constants and frequency threshold are uniform for $z_0\in\mathcal K$, $\theta_0\in\Torus$, $0\le\sigma\le\bar\sigma$, and $t\ge0$.
\end{theorem}
\begin{proof}
Proposition~\ref{prop:deterministic-loop} gives $\Phi_\omega(\theta_0,z_0)<r$ and \eqref{eq:det-Phi}. On $\overline{\mathcal D_R^\omega(\theta)}$, write $y=\Psi_\omega(\theta,z)$. Section~\ref{sec:centered} places $z$ in a fixed compact subset of $\mathcal O_1$ and $y$ in $\mathcal S_{R_0}$.

Using \eqref{eq:gamma-bounds}, we choose uniform bounds $\|D_z\Psi_\omega\|\le J_\Psi$ and $\|D_z^2\Psi_\omega\|\le H_\Psi$ . Let $G_V,H_V$ bound $|\nabla V|,\|D^2V\|$ on $\mathcal S_{R_0}$, and let $S_*=\max_{\overline{\mathcal O_1}}\|S_0\|_F$. The chain rule applied to \eqref{eq:Phi} gives
\[
 \nabla_z\Phi_\omega=(D_z\Psi_\omega)\transpose\nabla V(y)
\]
and
\[
\begin{split}
 D_z^2\Phi_\omega[v,w]={}D^2V(y)&[D_z\Psi_\omega v,D_z\Psi_\omega w]\\
 &+\nabla V(y)\transpose D_z^2\Psi_\omega[v,w].
\end{split}
\]
Thus $\|D_z^2\Phi_\omega\|\le H_\Phi:=J_\Psi^2H_V+G_VH_\Psi$. Equation~\eqref{eq:extended-generator} yields
\begin{equation}\label{eq:closed-generator}
\begin{split}
 \cA_{\omega,\sigma}\Phi_\omega
 &=\cA_{\omega,0}\Phi_\omega
       +\tfrac12\sigma^2\tr(S_0S_0\transpose D_z^2\Phi_\omega)\\
 &\le-\lambda\Phi_\omega+c_0\sigma^2,
\end{split}
\end{equation}
where $c_0=S_*^2H_\Phi/2$. The estimate $|\nabla V(y)|\le L_V|y-z^\star|$ and \eqref{eq:certificate-positive} similarly give
\begin{equation}\label{eq:closed-qv}
 \sigma^2|S_0\transpose\nabla_z\Phi_\omega|^2
 \le d_0\sigma^2\Phi_\omega
\end{equation}
where \(d_0=S_*^2J_\Psi^2L_V^2/c_1.\)
All bounds are independent of phase, $\sigma$, and sufficiently large $\omega$.

The compact containment and boundary value $\Phi_\omega=R$ allow us to apply Lemma~\ref{lem:stopped} to $(t,z)\mapsto\Phi_\omega(\theta_0+\omega t,z)$. This gives
\begin{equation}\label{eq:physical-killed-Phi}
\begin{split}
& \E[\Phi_\omega(\theta_t,Z_t^\omega)\ind_{\{t<\tau_R^\omega\}}] \\
&\quad \le e^{-\lambda t}\Phi_\omega(\theta_0,z_0)+\frac{c_0\sigma^2}{\lambda}(1-e^{-\lambda t}).
\end{split}
\end{equation}
The lower bound in \eqref{eq:orbital-comparability}
and \eqref{eq:physical-killed-Phi} give
\[
\begin{split}
 &\E\!\left[
 |Z_t^\omega-z_\omega^\circ(\theta_t)|^2
 \ind_{\{t<\tau_R^\omega\}}
 \right]\\
 &\quad\le
 \frac{e^{-\lambda t}}{a_1}\Phi_\omega(\theta_0,z_0)
 +\frac{c_0\sigma^2}{\lambda a_1}(1-e^{-\lambda t}).
\end{split}
\]
At the initial state, the upper bound gives
$\Phi_\omega(\theta_0,z_0)
\le a_2|z_0-z_\omega^\circ(\theta_0)|^2$.
Substitution proves \eqref{eq:orbital-moment}
with $M=a_2/a_1$ and $C\ge c_0/(\lambda a_1)$.

For $\sigma>0$, apply \eqref{eq:abstract-exit}
with $Q_0=\Phi_\omega(\theta_0,z_0)\le r$.
Since $R-Q_0\ge R-r$, replacing $Q_0$ by $r$
in that upper bound proves \eqref{eq:main-exit}.

Finally, Lemma~\ref{lem:orbit} gives
$\sup_\theta|z_\omega^\circ(\theta)-z^\star|
\le C_o\omega^{-1/2}$. The inequality
$|x+y|^2\le2|x|^2+2|y|^2$ therefore yields
\[
 |Z_t^\omega-z^\star|^2
 \le
 2|Z_t^\omega-z_\omega^\circ(\theta_t)|^2
 +\frac{2C_o^2}{\omega}
\]
before exit, and
\[
 |z_0-z_\omega^\circ(\theta_0)|^2
 \le 2|z_0-z^\star|^2+\frac{2C_o^2}{\omega}.
\]
Taking the pre-exit expectation in the first inequality
and substituting the orbital estimate and the second
inequality gives
\[
\begin{split}
 &\E\!\left[
 |Z_t^\omega-z^\star|^2
 \ind_{\{t<\tau_R^\omega\}}
 \right]\\
 &\quad\le
 4Me^{-\lambda t}|z_0-z^\star|^2
 +\frac{2c_0}{\lambda a_1}\sigma^2
 +\frac{(4M+2)C_o^2}{\omega}.
\end{split}
\]
We then choose
\[
 C\ge\max\left\{
 4M,\,
 \frac{2c_0}{\lambda a_1},\,
 (4M+2)C_o^2
 \right\}.
\]
The preceding right-hand side is then bounded by
\[
 Ce^{-\lambda t}|z_0-z^\star|^2
 +C(\sigma^2+\omega^{-1}),
\]
which establishes \eqref{eq:physical-moment}.
This choice of $C$ is independent of time, the
initial state and phase, and the frequency and
noise intensity within the stated ranges.
\end{proof}

Since $\nabla h(x^\star)=0$, a Taylor expansion near $x^\star$ and the compactness away from it give $|h(x)-h^\star|\le C_h|x-x^\star|^2$ on the physical analysis region. The expected absolute cost error before exit therefore satisfies \eqref{eq:physical-moment} with a changed constant.

The next corollary strengthens the regulation estimate
to a trajectory-level statement. With high probability,
the closed-loop state enters a prescribed neighborhood
of $z^\star$ by a fixed capture time and remains there
throughout a subsequent interval. The proof first bounds
the noise fluctuations that could prevent capture, then
applies the exit estimate of Lemma~\ref{lem:stopped}
from the state reached at the capture time. Both failure
probabilities decrease exponentially with the inverse
noise variance, allowing exponentially long residence
intervals.

\begin{corollary}\label{cor:capture}
Under Theorem~\ref{thm:main}, for each $\rho>0$ there exist $T_a>0$, $\omega_\rho\ge\omega_{\mathcal K}$, and $0<\sigma_\rho\le\bar\sigma$ such that
\begin{equation}\label{eq:capture}
\begin{split}
 \Prob\Big\{&\tau_R^\omega>T_a+T,\
 \sup_{T_a\le t\le T_a+T}|Z_t^\omega-z^\star|<\rho\Big\}\\
 &\ge1-e^{-b_{\rm cap}/\sigma^2}-e^{-b_1/\sigma^2}
                      -C_0T e^{-b_2/\sigma^2}.
\end{split}
\end{equation}
The constants $b_{\rm cap},b_1,b_2,C_0>0$ are uniform for $z_0\in\mathcal K$, $\theta_0\in\Torus$, $\omega\ge\omega_\rho$, $0<\sigma\le\sigma_\rho$, and finite $T\ge0$. For $T=T_*e^{a/\sigma^2}$, with fixed $T_*>0$ and $0<a<b_2$, the failure bound is at most $(2+C_0T_*)e^{-b_a/\sigma^2}$, where $b_a=\min\{b_{\rm cap},b_1,b_2-a\}>0$.
\end{corollary}
\begin{proof}
We start by choosing $0<r_a<R_a<R$ with
$\sqrt{R_a/a_1}<\rho/2$. By Lemma~\ref{lem:orbit},
we can increase $\omega_\rho$ such that
$\sup_\theta|z_\omega^\circ(\theta)-z^\star|<\rho/2$.
The lower bound in \eqref{eq:orbital-comparability}
then ensures that
$\mathcal D_{R_a}^\omega(\theta)$ lies inside the
$\rho$-ball about $z^\star$ for every phase.
We then choose $T_a>0$ such that
$e^{-\lambda T_a}r<r_a/4$.

Let $Q_t=\Phi_\omega(\theta_t,Z_t^\omega)$ up
to the outer exit time $\tau_R^\omega$. We first
estimate the probability of the capture event
$\mathcal G:=\{\tau_R^\omega>T_a,\ Q_{T_a}<r_a\}$.
To control the accumulated noise before $T_a$,
we define
\[
 M_t
 :=
 \sigma\int_0^{t\wedge\tau_R^\omega}
 e^{\lambda s}
 \nabla_z\Phi_\omega(\theta_s,Z_s^\omega)\transpose
 S_0(Z_s^\omega)\,\dd B_s.
\]
Applying It\^o's formula to $e^{\lambda t}Q_t$
up to $\tau_R^\omega$ and using
\eqref{eq:closed-generator} gives
\begin{equation}\label{eq:capture-pathwise}
 Q_{t_R}
 \le e^{-\lambda t_R}Q_0
 +\frac{c_0\sigma^2}{\lambda}
 +e^{-\lambda t_R}M_t
\end{equation}
for \(t_R:=t\wedge\tau_R^\omega.\)
Since $Q_s<R$ before exit, the expression \eqref{eq:closed-qv}
gives
\[
 \langle M\rangle_{T_a}
 \le d_0R\sigma^2\int_0^{T_a}e^{2\lambda s}\,\dd s
 \le H_a\sigma^2,
\]
where
$H_a:=\max\{1,d_0R(e^{2\lambda T_a}-1)/(2\lambda)\}$.

We then set $\Delta_a=\tfrac14\min\{R-r,r_a\}$ and
choose $0<\sigma_\rho\le\bar\sigma$ such that
$c_0\sigma_\rho^2/\lambda\le\Delta_a$.
On the event $\{\sup_{t\le T_a}M_t<\Delta_a\}$,
inequality \eqref{eq:capture-pathwise} implies
$Q_{t_R}\le r+2\Delta_a<R$ for every $t\le T_a$.
An outer exit would require
$Q_{\tau_R^\omega}=R$, so $\tau_R^\omega>T_a$.
At the capture time, the same estimate gives
\[
 Q_{T_a}
 \le e^{-\lambda T_a}r+2\Delta_a
 <\frac34r_a<r_a.
\]
Thus a failure of $\mathcal G$ requires the
martingale $M_t$ to reach $\Delta_a$ before $T_a$.

For $\gamma>0$, the process
$\exp(\gamma M_t-\gamma^2\langle M\rangle_t/2)$
is a nonnegative supermartingale
\cite{Howard2020}. It is stopped when $M_t$ first
reaches $\Delta_a$, or at $T_a$. We then choose
$\gamma=\Delta_a/(H_a\sigma^2)$. This gives
\[
 \Prob(\mathcal G^c)
 \le
 \Prob\!\left\{\sup_{t\le T_a}M_t\ge\Delta_a\right\}
 \le
 \exp\!\left(-\frac{\Delta_a^2}{2H_a\sigma^2}\right).
\]
Hence the capture failure probability is at most
$e^{-b_{\rm cap}/\sigma^2}$, with
$b_{\rm cap}:=\Delta_a^2/(2H_a)$.

It remains to monitor the departures after capture.
On $\mathcal G$, let $\tau_a^\omega$ be the first
exit from $\mathcal D_{R_a}^\omega(\theta_t)$
after $T_a$. The bounds
\eqref{eq:closed-generator}--\eqref{eq:closed-qv}
hold on this smaller sublevel, uniformly in phase.
Applying Lemma~\ref{lem:stopped} conditionally on
$\mathcal F_{T_a}$, with $Q_{T_a}<r_a$, yields
\[
\begin{split}
 &\Prob\!\left\{
 \tau_a^\omega\le T_a+T
 \,\middle|\,\mathcal F_{T_a}
 \right\}\\
 &\quad\le
 e^{-\eta(R_a-r_a)/\sigma^2}
 +C_0T e^{-\eta R_a/\sigma^2}
 \qquad\text{on }\mathcal G.
\end{split}
\]
On $\mathcal G\cap\{\tau_a^\omega>T_a+T\}$,
the trajectory remains in the $\rho$-ball during
$[T_a,T_a+T]$ and does not leave the outer sublevel
before the interval ends. Adding the capture
failure bound to the conditional exit bound proves
\eqref{eq:capture}, with
$b_1:=\eta(R_a-r_a)$ and $b_2:=\eta R_a$.

Finally, for $T=T_*e^{a/\sigma^2}$ and $0<a<b_2$,
the failure bound becomes
\[
\begin{split}
 &e^{-b_{\rm cap}/\sigma^2}
 +e^{-b_1/\sigma^2}
 +C_0T_*e^{-(b_2-a)/\sigma^2}\\
 &\quad\le
 (2+C_0T_*)e^{-b_a/\sigma^2},
\end{split}
\]
where $b_a=\min\{b_{\rm cap},b_1,b_2-a\}>0$.
All choices use the uniform bounds of
Theorem~\ref{thm:main}, so the capture time and
constants are independent of the initial state,
initial phase, and frequency above $\omega_\rho$.
\end{proof}
\subsection{Time-uniform approximation under incremental stability}
The estimate \eqref{eq:finite-window} bounds an expected trajectory supremum on each fixed interval. We now bound the mean-square discrepancy at every time with one constant. This requires decay of the separation between averaged trajectories driven by the same Brownian motion.

\begin{assumption}\label{ass:incremental}
For the fixed gains, region $\mathcal O_2$, and noise bound $\bar\sigma$, there exist $P\in\Spp^{n+m}$ and $\alpha>0$ such that, for $z,z'\in\mathcal O_2$ and $0\le\sigma\le\bar\sigma$,
\begin{equation}\label{eq:incremental}
\begin{split}
 &2(z-z')\transpose P[F_T(z)-F_T(z')]\\
 &\quad+\sigma^2\|P^{1/2}[S_0(z)-S_0(z')]\|_F^2
 \le-\alpha|z-z'|_P^2,
\end{split}
\end{equation}
where $|e|_P^2=e\transpose Pe$.
\end{assumption}
\begin{remark}
This condition is stronger than the equilibrium certificate on the same region. Its verification requires structural information or uncertainty bounds that are not supplied by the scalar cost alone. It is used only for the approximation result. On convex $\mathcal O_2$, a sufficient test is
\[
 DF_T\transpose P+PDF_T+
 \bar\sigma^2\sum_{j=1}^{d_B}DS_{0,j}\transpose PDS_{0,j}\preceq-\alpha P.
\]
Here $S_{0,j}$ is the $j$th diffusion column. Integration along a segment and Jensen's inequality imply \eqref{eq:incremental}. A Hurwitz target linearization gives this test on a sufficiently small neighborhood for small enough $\bar\sigma$, not automatically on the prescribed region.
\end{remark}

\begin{theorem}\label{thm:uniform}
Let Assumptions~\ref{ass:smooth}--\ref{ass:certificate} and \ref{ass:incremental} hold. Fix an open $\mathcal U$ with $\mathcal K\subset\mathcal U\Subset\mathcal O_1$, and couple \eqref{eq:complete} and \eqref{eq:average-sde} with the same initial state and Brownian motion. Put $\tau_c^\omega=\tau_{\mathcal U}(Z^\omega)\wedge\tau_{\mathcal U}(\bar Z)$. There exist $\omega_0,C>0$ such that
\begin{equation}\label{eq:uniform-approximation}
 \sup_{\substack{z_0\in\mathcal K,\,\theta_0\in\Torus\\t\ge0}}
 \E[|Z_t^\omega-\bar Z_t|^2\ind_{\{t<\tau_c^\omega\}}]\le C/\omega
\end{equation}
for $\omega\ge\omega_0$ and $0\le\sigma\le\bar\sigma$.
\end{theorem}
\begin{proof}
Before $\tau_c^\omega$, define $\zeta_t=\Gamma_\eps^{-1}(\theta_t,Z_t^\omega)$ and $e_t=\zeta_t-\bar Z_t$. The buffer and \eqref{eq:gamma-bounds} keep both arguments in $\mathcal O_2$ for small $\eps$. Their error equation is
\[
 \dd e_t=[F_T(\zeta_t)-F_T(\bar Z_t)+r_t]\dd t
                       +\sigma(D_t+J_t)\dd B_t,
\]
where $r_t=\eps R_\eps(\theta_t,\zeta_t)+\sigma^2B_\eps(\theta_t,\zeta_t)$, $D_t=S_0(\zeta_t)-S_0(\bar Z_t)$, and $J_t=\widehat S_\eps(\theta_t,\zeta_t)-S_0(\zeta_t)$. Lemma~\ref{lem:normal} gives $|r_t|+\|J_t\|_F\le C\eps$ for $\sigma\le\bar\sigma$. A Lipschitz bound for the extension of $S_0$ gives $\|D_t\|_F\le L_S|e_t|$.

For $W_t=e_t\transpose Pe_t$, the generator of the coupled process satisfies
\[
\begin{split}
 \cL_{\rm pair}W
 &=2e_t\transpose P[F_T(\zeta_t)-F_T(\bar Z_t)+r_t]\\
 &\quad+\sigma^2\|P^{1/2}(D_t+J_t)\|_F^2\\
 &\le-\alpha W_t+2e_t\transpose Pr_t\\
 &\quad+2\sigma^2\langle P^{1/2}D_t,P^{1/2}J_t\rangle_F
       +\sigma^2\|P^{1/2}J_t\|_F^2,
\end{split}
\]
where \eqref{eq:incremental} bounds the terms containing only $F_T,D_t$, and $\langle A,B\rangle_F=\tr(A\transpose B)$. The remainder is at most $C\eps\sqrt{W_t}+C\eps^2$. Young's inequality yields
\begin{equation}\label{eq:joint-drift}
 \cL_{\rm pair}W\le-\tfrac12\alpha W_t+C_1\eps^2,
 \qquad t<\tau_c^\omega.
\end{equation}
The stopped stochastic integral is square integrable on each finite interval. Applying \eqref{eq:joint-drift} to $e^{\alpha(t\wedge\tau_c^\omega)/2}W_{t\wedge\tau_c^\omega}$ and dropping its nonnegative exit contribution gives
\[
 \E[W_t\ind_{\{t<\tau_c^\omega\}}]
 \le e^{-\alpha t/2}W_0+\frac{2C_1}{\alpha}\eps^2(1-e^{-\alpha t/2}).
\]
Let $C_\Gamma$ denote a constant in \eqref{eq:gamma-bounds}. The common initialization implies $W_0\le\|P\|C_\Gamma^2\eps^2$ by \eqref{eq:gamma-bounds}. The same bound and $Z_t^\omega=\Gamma_\eps(\theta_t,\zeta_t)$ give
\[
 |Z_t^\omega-\bar Z_t|^2\le2C_\Gamma^2\eps^2+2W_t/\lambda_{\min}(P).
\]
Combining these inequalities proves \eqref{eq:uniform-approximation}, uniformly in time, phase, initial state, and $\sigma\le\bar\sigma$.
\end{proof}

To relate this pre-exit estimate to a growing interval, we choose $\mathcal S_{R_0}\subset\mathcal U\Subset\mathcal O_1$. For a sufficiently large frequency, both $\overline{\mathcal D_R^\omega(\theta)}$ and $\mathcal S_R=\{V\le R\}$ lie in $\mathcal U$. With $\bar\tau_R=\inf\{t:V(\bar Z_t)\ge R\}$, the union bound gives
\[
 \Prob\{\tau_c^\omega\le T\}
 \le\Prob\{\tau_R^\omega\le T\}+\Prob\{\bar\tau_R\le T\}.
\]
Apply \eqref{eq:main-exit} to the first term and Lemma~\ref{lem:stopped} with \eqref{eq:average-generator-bound} to the second. Taking the smaller exponential coefficient gives
\[
 \Prob\{\tau_c^\omega\le T\}\le\delta_\sigma(T)
 :=2e^{-\eta_*(R-r)/\sigma^2}+C_*T e^{-\eta_*R/\sigma^2}
\]
for constants $\eta_*,C_*>0$. Since $\{\tau_c^\omega>T\}\subset\{t<\tau_c^\omega\}$ for $t\le T$, \eqref{eq:uniform-approximation} implies
\[
 \sup_{0\le t\le T}\E[|Z_t^\omega-\bar Z_t|^2\mid\tau_c^\omega>T]
 \le\frac{C}{\omega[1-\delta_\sigma(T)]}
\]
whenever $\delta_\sigma(T)<1$. For $T=T_*e^{a/\sigma^2}$ with $0<a<\eta_*R$, $\delta_\sigma(T)\to0$ as $\sigma\downarrow0$. The supremum here is outside the expectation, as in Theorem~\ref{thm:uniform}.
The estimate bounds the conditional mean-square error  at every time in $[0,T]$, given that both trajectories remain in $\mathcal U$ throughout that interval. The supremum is taken over these expectations, rather than over the error along each individual trajectory.
\section{Brownian-Particle Confinement in a Hydrodynamic Saddle Trap}
\label{secVIII}
\label{sec:trap}

\subsection{Physical model and measured objective}
Consider a planar hydrodynamic trap with controlled outlet-axis stagnation-point position $x_{\rm s}=x_{\rm s0}+\beta_c u$ \cite{TrapTanyeri2010,TrapTanyeri2011,TrapTanyeri2013}. The physical particle position $(x_{\rm p},y_{\rm p})$ satisfies the local overdamped model
\[
\begin{aligned}
 \dd x_{\rm p}&=\dot\epsilon_f(x_{\rm p}-x_{\rm s})\,\dd\tau
                 +\sqrt{2D_B}\,\dd W_{1,\tau},\\
 \dd y_{\rm p}&=-\dot\epsilon_f y_{\rm p}\,\dd\tau
                 +\sqrt{2D_B}\,\dd W_{2,\tau},
\end{aligned}
\]
where $\dot\epsilon_f>0$ is the strain rate, $D_B>0$ is the translational diffusivity, and $W_1,W_2$ are independent Wiener processes. The offset $x_{\rm s0}$ is unknown and $\beta_c\ne0$ is the signed command sensitivity. All coefficients are fixed during each realization.

Normalize position by the detection-spot width $w>0$ and time by a reference strain rate $\dot\epsilon_0>0$. With $x=x_{\rm p}/w$, $y=y_{\rm p}/w$, and $t=\dot\epsilon_0\tau$, the dynamics become
\begin{equation}\label{eq:trap-model}
\begin{aligned}
 \dd x_t&=(ax_t+d+bu_t)\,\dd t+\sigma\,\dd B_{1,t},\\
 \dd y_t&=-ay_t\,\dd t+\sigma\,\dd B_{2,t},
\end{aligned}
\end{equation}
where $a=\dot\epsilon_f/\dot\epsilon_0$, $d=-ax_{\rm s0}/w$, $b=-a\beta_c/w$, and $\sigma^2=2D_B/(\dot\epsilon_0w^2)$. The controller receives only the optical score
\begin{equation}\label{eq:trap-cost}
 Y_t=h(x_t,y_t)
 =h_0-I_0\exp\!\left[-\frac{x_t^2+y_t^2}{2}\right],
\end{equation}
with constant background $h_0$ and peak intensity $I_0>0$. Minimizing this score maximizes detected intensity, and the controller forms no position estimate. This idealized instantaneous measurement omits the photon-arrival sensing and spatial scanning used for position estimation in single-molecule ABEL traps \cite{TrapWang2011}.

The plant has state $X=(x,y)^\top\in\R^2$, input vector $G=(b,0)^\top$, and normalized diffusion matrix $\Sigma_0=I_2$. Thus, noise acts along the non-actuated, passively stable $y$ direction. The feasible nominal minimum is $X^\star=0$, with sustaining input $u^\star=-d/b$, value $h^\star=h_0-I_0$, and Hessian $\nabla^2h(0)=I_0I_2\succ0$. Writing $\phiT(x,y)=\exp[-(x^2+y^2)/2]$, the input sensitivity is $\ell(x,y)=bI_0x\phiT(x,y)$. The scripts denote the strain-rate ratio $a$ by \texttt{eps}, whereas the averaging parameter is $\omega^{-1/2}$.

\subsection{Controller and integral compensation}
We use \eqref{eq:controller} with $\kappa_1=1$, $k_P=4$, and $k_I=1$. With $\kappa_P=k_Pb^2I_0$, the averaged diffusion is
\begin{equation}\label{eq:trap-average}
\begin{aligned}
 \dd\bar x_t&=[a\bar x_t+d+b\bar u_{b,t}
       -\kappa_P\bar x_t\phiT(\bar x_t,\bar y_t)]\,\dd t
       +\sigma\,\dd B_{1,t},\\
 \dd\bar y_t&=-a\bar y_t\,\dd t+\sigma\,\dd B_{2,t},\\
 \dd\bar u_{b,t}&=-k_IbI_0\bar x_t\phiT(\bar x_t,\bar y_t)\,\dd t.
\end{aligned}
\end{equation}

To isolate integral compensation, we compare against the matched proportional-only controller obtained by fixing $u_b=0$ while retaining the same proportional gain and probing signal. At $y=0$, its deterministic average is
\begin{equation}\label{eq:trap-ponly}
 \dot x=ax+d-\kappa_Pxe^{-x^2/2}.
\end{equation}
The drift at the optical minimum equals $d$, so this point is not an equilibrium when $d\ne0$. For $\kappa_P>a$ and sufficiently small $d$, the locally stable equilibrium satisfies $x_d=d/(\kappa_P-a)+O(d^3)$. The dual-mode deterministic average admits the optical minimum as an equilibrium through the sustaining bias $u_b^\star=-d/b$.

\subsection{Regional target stability}
Consider the uncertainty family $a\in[0.8,1.2]$, $|b|\in[0.9,1.1]$, $I_0\in[0.9,1.1]$, and $|d|\le0.25$. Define $\varrho=b(u_b-u^\star)=bu_b+d$, $\xi=\col(x,y,\varrho)$, and $k=b^2I_0=\kappa_I$. Appendix~\ref{app:trap-certificate} shows that
\begin{equation}\label{eq:trap-V}
 V(\xi)=\tfrac12(x^2+y^2)+\frac{\varrho^2}{2k}
                -\tfrac14x\varrho
\end{equation}
satisfies the uniform bounds
\begin{equation}\label{eq:trap-certified}
 0.298|\xi|^2\le V(\xi)\le0.749|\xi|^2
\end{equation}
and $\dot V\le-0.1549|\xi|^2$ along the noise-free target for $|\xi|<0.7$. In the original error coordinates $e=\col(x,y,u_b-u^\star)$, the corresponding constants in Definition~\ref{def:certificate} are $c_1=0.24138$, $c_2=0.90629$, and $c_3=0.125469$.

The outer sublevel $\{V\le R_0\}$, with $R_0=0.14$, lies strictly inside $|\xi|<0.7$, since $R_0<0.298(0.7)^2$. At the common initialization $(x_0,y_0,u_{b,0})=(0.2,0.1,0)$, the nominal choice $k=1$, $d=0.25$ gives $V(\xi_0)=0.04375$. Over the entire uncertainty family, $V(\xi_0)\le0.080367$, so $r_0=0.09<R_0$ covers this initialization uniformly.

Since the noise acts additively on $(x,y)$ and $V_{xx}+V_{yy}=2$, the averaged generator and quadratic-variation rate satisfy
\begin{equation}\label{eq:trap-generator}
\begin{aligned}
 \overline{\mathcal L}_\sigma V
 &\le-\frac{0.1549}{0.749}V+\sigma^2,\\
 q_\sigma(V)
 &=\sigma^2[(x-\varrho/4)^2+y^2]
 \le2\sigma^2V
\end{aligned}
\end{equation}
on the certified sublevel. These estimates verify the regional target certificate for Theorem~\ref{thm:main} and the averaged bounds used in Lemma~\ref{lem:stopped}.

\subsection{Simulation protocol}
Simulations use Euler--Maruyama integration with nominal parameters $a=b=I_0=1$, $d=0.25$, and $h_0=\theta_0=0$, and the initialization specified above. Table~\ref{tab:trap-protocol} summarizes the ensemble sizes and simulation horizons. Paired dual-mode and proportional-only runs use identical Wiener increments.

\begin{table}[t]
\centering
\caption{Simulation settings. $N$ is the number of paths and $T$ is the normalized time horizon.}
\label{tab:trap-protocol}
\begin{tabular}{@{}lrr@{}}
\toprule
Experiment & $N$ & $T$\\
\midrule
Nominal, offset, actuator sign & 128 & 20\\
Frequency sweep, $\sigma=0$ & 1 & 20\\
Frequency sweep, $\sigma>0$ & 32 & 20\\
Large-offset stress test & 64 & 12\\
Strip-exit study & 1024 & 200\\
\bottomrule
\end{tabular}
\end{table}

Curves labelled ``period averaged'' are moving averages of the simulated paths, rather than trajectories of the averaged diffusion \eqref{eq:trap-average}. Shaded bands span the pointwise empirical fifth to ninety-fifth percentiles. Proportional-only equilibrium references are computed from \eqref{eq:trap-ponly}.

\subsection{Regulation, offset compensation, and frequency dependence}
\begin{figure}[t]
\centering
\includegraphics[width=\columnwidth]{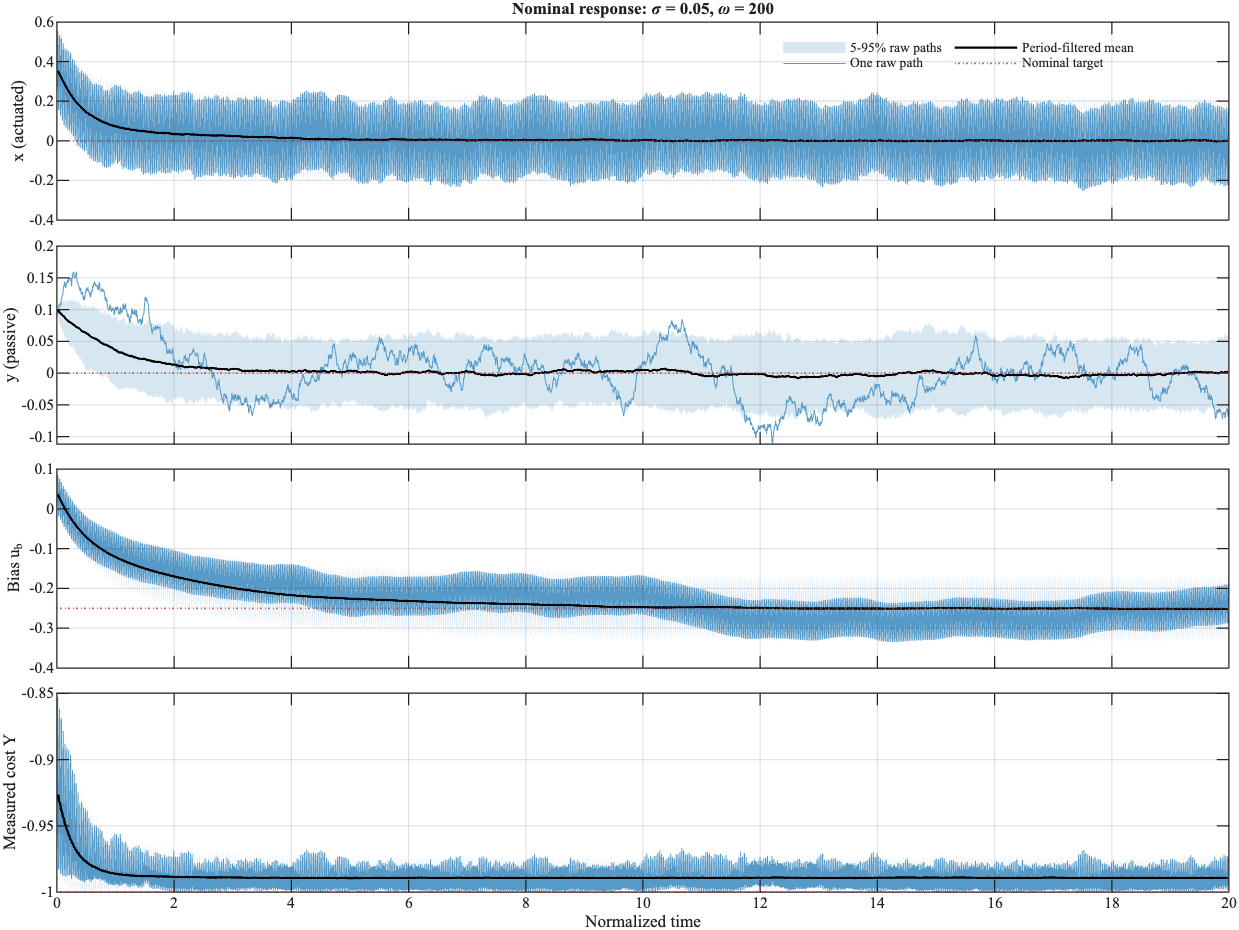}
\caption{Nominal trap response over 128 paths, showing positions, bias, and cost. Thin curves show one realization and thick curves show period-filtered ensemble means. Shading spans pointwise 5--95\% empirical quantiles of the unfiltered paths. Dotted lines mark nominal targets.}
\label{fig:trap-A}
\end{figure}

For $\sigma=0.05$ and $\omega=200$, the filtered mean positions settle near the spot center and the mean bias approaches $u^\star=-0.25$ (Fig.~\ref{fig:trap-A}). The cost remains near its minimum, with residual fluctuations from probing and Brownian motion.

\begin{figure*}[t]
\centering
\includegraphics[width=0.94\textwidth]{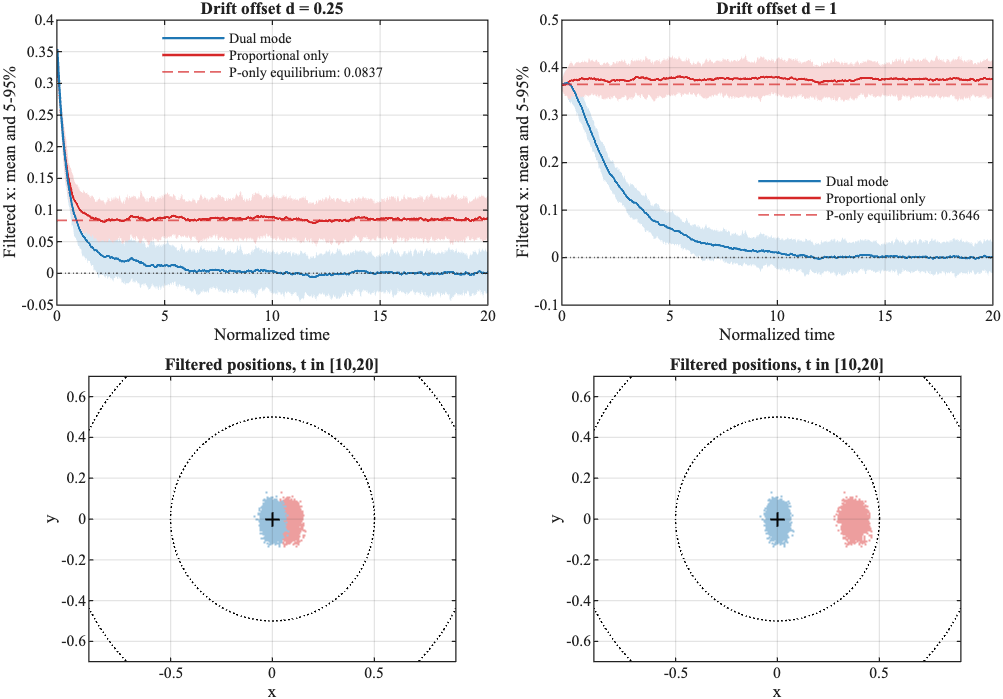}
\caption{Dual-mode and matched frozen-bias simulations for $d=0.25$ and $d=1$. Upper panels: filtered ensemble mean and empirical 5--95\% bands, with the stable proportional-only equilibrium indicated. Lower panels: filtered particle positions over $t\in[10,20]$. Circles mark normalized radii $0.5$ and $1$ for visualization, not a stopping criterion used in these runs. The $d=1$ case is outside the certified offset family.}
\label{fig:trap-B}
\end{figure*}

In Fig.~\ref{fig:trap-B}, the filtered proportional-only means agree with the equilibria predicted by \eqref{eq:trap-ponly}, approximately $0.0837$ for $d=0.25$ and $0.3646$ for $d=1$. The dual-mode means remain near zero. The $d=1$ case is a stress test outside the certified uncertainty family.

\begin{figure*}[t]
\centering
\includegraphics[width=\textwidth]{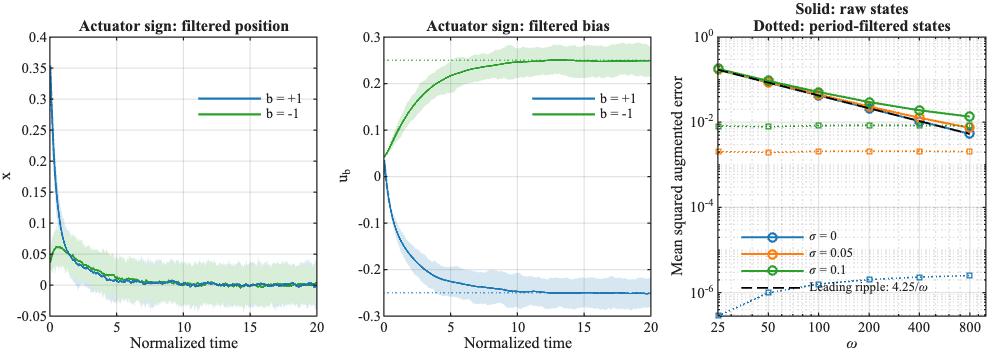}
\caption{Left and center: one-period-filtered position and bias for both actuator signs. Right: reported sample/time mean squared augmented error over $t\in[10,20]$. Solid curves use the stored raw state and dotted curves use its moving average. The dashed line is the leading ripple energy \eqref{eq:trap-ripple}. The frequency sweep uses one noise-free path or 32 noisy paths per point, not 128.}
\label{fig:trap-D}
\end{figure*}

Figure~\ref{fig:trap-D} shows regulation for either actuator sign, with the bias approaching $+0.25$ when $b=-1$. For nominal $b=1$, the frequency sweep uses $\omega\in\{25,50,100,200,400,800\}$ and $\sigma\in\{0,0.05,0.1\}$. The reported error averages $|Z-z^\star|^2$ over paths and $t\in[10,20]$, where $Z=\col(x,y,u_b)$ and $z^\star=\col(0,0,-0.25)$.

The leading-order ripple energy is
\begin{equation}\label{eq:trap-ripple}
 E_{\rm ripple}(\omega)
 =\frac{k_Pb^2+k_I^2/k_P}{\omega}
 =\frac{4.25}{\omega}.
\end{equation}
The nominal linearized averaged diffusion has stationary covariance $P_{\rm lin}=\sigma^2\diag(1/6,1/2,1/6)$. Adding its trace gives the local reference $4.25/\omega+5\sigma^2/6$, which agrees with the simulated errors at $\omega=200$ (Table~\ref{tab:trap-mse}) and supports the ripple--noise interpretation.

\begin{table}[t]
\centering
\caption{Raw mean squared error at $\omega=200$ and the local analytical reference.}
\label{tab:trap-mse}
\begin{tabular}{@{}rrr@{}}
\toprule
$\sigma$ & Sample/time mean & $4.25/\omega+5\sigma^2/6$\\\midrule
$0$    & $0.02126$ & $0.02125$\\
$0.05$ & $0.02332$ & $0.02333$\\
$0.10$ & $0.02950$ & $0.02958$\\
\bottomrule
\end{tabular}
\end{table}

\subsection{Large-offset capture}

\begin{figure*}[t]
\centering
\includegraphics[scale=0.75]{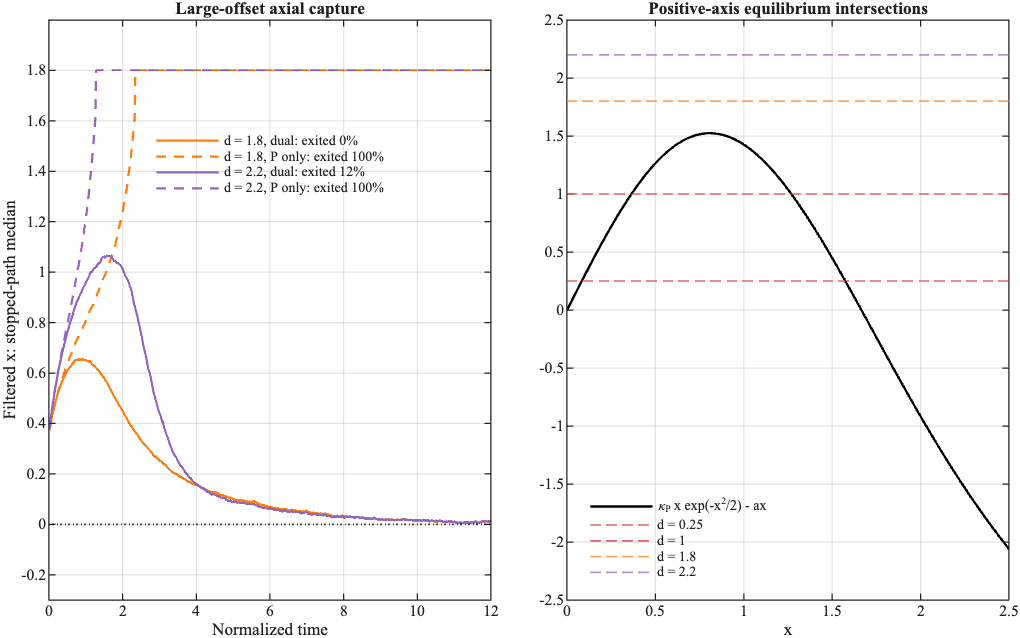}
\caption{Large-offset stress tests, 64 paths per case. Left: medians of the one-period-filtered, grid-stopped positions. ``Lost'' denotes first exit from $|x|<1.8$ by $T=12$ and is rounded to a whole percentage. Stopped paths remain in the plotted medians. Right: positive-axis intersections for the proportional-only average, illustrating disappearance of its stable branch. This panel does not display every equilibrium on $\R$.}
\label{fig:trap-E}
\end{figure*}

The offsets $d=1.8$ and $2.2$ lie outside the certified uncertainty family. In 64-path runs to $T=12$, all proportional-only paths exit the strip $|x|<1.8$ at both offsets. The dual-mode exit fractions are $0\%$ and approximately $19\%$, respectively (Fig.~\ref{fig:trap-E}), illustrating improved transient capture under large offsets.

Paths are frozen at their first recorded exit. The plotted medians include these paths and moving-average post-processing, so they describe the full stopped ensemble rather than surviving paths alone. The right panel shows only $x\ge0$ and omits the unstable equilibria at $x<0$.

The simulations assume ideal stagnation-point actuation and instantaneous intensity feedback. At the nominal settings, the probing signal requires a stagnation-point displacement amplitude of $40w$ at frequency $\omega\dot\epsilon_0/(2\pi)$. A physical implementation must accommodate this excursion and probing frequency.
\section{Conclusion}\label{secIX}
A dual-mode extremum-seeking controller with cosine modulation was analyzed for unknown input-affine controlled diffusions. An exact deterministic periodic response separates probing ripple from process noise in the Lyapunov analysis. Regional moment bounds and exponential capture and residence estimates establish practical regulation from certified initial sets over horizons that grow as the noise decreases. Under an additional incremental condition, the same-noise approximation is uniform in mean square before common exit. The hydrodynamic trap example verifies a regional target certificate and illustrates regulation and compensation of unknown offsets using scalar optical feedback. Future work will examine the certified frequency requirement under sensing and actuation constraints.

\appendices
\section{Verification of the Regional Trap Certificate}
\label{app:trap-certificate}

Let $k=b^2I_0$ and $\phiT=e^{-(x^2+y^2)/2}$. Differentiating \eqref{eq:trap-V} along the noise-free dynamics in \eqref{eq:trap-average} gives
\begin{equation}\label{eq:trap-Vdot}
 \dot V=-ay^2-
 \begin{bmatrix}x&\varrho\end{bmatrix}
 M(a,k,\phiT)
 \begin{bmatrix}x\\\varrho\end{bmatrix},
\end{equation}
where $c_\times=1-a/4+(k-1)\phiT$ and
\begin{equation}\label{eq:trap-M}
 M(a,k,\phiT)=
 \begin{bmatrix}
 \tfrac{15}{4}k\phiT-a&-c_\times/2\\
 -c_\times/2&\tfrac14
 \end{bmatrix}.
\end{equation}

On $|\xi|<0.7$, $\phiT\ge e^{-0.245}>0.7827$. The matrix $M$ is affine in each variable over the box $a\in[0.8,1.2]$, $k\in[0.729,1.331]$, and $\phiT\in[0.7827,1]$, so every value is a convex combination of its eight vertex matrices. At each vertex, exact rational evaluation gives
$M_{11}-0.1549>0.7848$ and
$\det(M-0.1549I_2)>3.02\times10^{-4}$.
Hence $M\succeq0.1549I_2$ throughout the box. Since $a\ge0.8$, \eqref{eq:trap-Vdot} yields
\[
 \dot V\le-0.1549|\xi|^2
\]
throughout $|\xi|<0.7$, proving the required decay bound.

Writing $V(\xi)=\xi^\top P(k)\xi$ gives
\[
 P(k)=\frac12
 \begin{bmatrix}
 1&0&-\tfrac14\\
 0&1&0\\
 -\tfrac14&0&1/k
 \end{bmatrix}.
\]
Since $P$ is affine in $1/k$, endpoint checks using Sylvester's criterion give $0.298I_3\prec P(k)\prec0.749I_3$, proving \eqref{eq:trap-certified}. Consequently, $\{V\le0.14\}\subset\{|\xi|<0.7\}$.

At the common initialization, $V_0=0.025+d^2/(2k)-0.05d$ is maximized at $d=-0.25$ and $k=0.729$, giving $V_0<0.080367<r_0=0.09$. For $e=\col(x,y,u_b-u^\star)$, the transformation $\xi=\diag(1,1,b)e$ gives $0.81|e|^2\le|\xi|^2\le1.21|e|^2$ and hence the certificate constants stated in the main text.

The additive-noise It\^o correction is $\sigma^2$. Moreover,
\[
 2V-[(x-\varrho/4)^2+y^2]
 =\left(\frac1k-\frac1{16}\right)\varrho^2\ge0,
\]
so $q_\sigma(V)\le2\sigma^2V$. Together with the decay estimate, these bounds establish \eqref{eq:trap-generator} and complete the regional certificate for Theorem~\ref{thm:main}.
\hfill$\square$

\section{Passive-Axis Bound on Radial Residence}
\label{app:trap-passive}

For $L>0$, define $\tau_{\rm spot}(L)=\inf\{t\ge0:x_t^2+y_t^2\ge L^2\}$ and $\tau_y(L)=\inf\{t\ge0:|y_t|\ge L\}$. These satisfy $\tau_{\rm spot}(L)\le\tau_y(L)$ pathwise.

For $\sigma>0$, the mean exit time $m(y)=\E_y\tau_y(L)$ of the passive coordinate solves $\frac{\sigma^2}{2}m''-aym'=-1$ with $m(\pm L)=0$. Integration gives
\[
 \E\tau_{\rm spot}(L)
 \le \E\tau_y(L)
 =\frac{\sqrt\pi}{a}
 \int_{\sqrt a|y_0|/\sigma}^{\sqrt aL/\sigma}
 e^{s^2}\operatorname{erf}(s)\,\dd s
\]
for an initial state inside the disk. This bound is independent of the controller. An optical intensity threshold $I_{\min}\in(0,I_0)$ corresponds to $L=\sqrt{2\log(I_0/I_{\min})}$.

\FloatBarrier
\bibliographystyle{IEEEtran}
\bibliography{references}

\end{document}